\documentclass[11pt,twoside,english]{article}
\usepackage{amssymb,amsmath}
\usepackage{amssymb,amsmath,babel}

\newtheorem{Theo_intro}{Theorem}
\newtheorem{Corollaire_intro}[Theo_intro]{Corollary}
\newtheorem{Lemme}{Lemma}[section]
\newtheorem{Theorem}[Lemme]{Theorem}

\newtheorem{Proposition}[Lemme]{Proposition}
\newtheorem{Corollaire}[Lemme]{Corollary}
\newtheorem{Definition}[Lemme]{Definition}

\newcommand{\NN}{\mathbb N}
\newcommand{\PP}{\mathbb P}
\newcommand{\ZZ}{\mathbb Z}
\newcommand{\QQ}{\mathbb Q}

\newcommand{\CC}{\mathbb C}

\newcommand{\FF}{\mathbb F}
\newcommand{\ie}{\emph{i.e.}}

\newcommand{\resp}{\emph{resp.}}
\newcommand{\sqm}[4]{\displaystyle{
\left({#1 \atop #3}{#2 \atop #4}\right)}}
\newcommand{\binomial}[2] {{\binom{#1}{#2}}}
\newcommand\CVD{{\hfill\hfil{\lower 2 pt\hbox{\vrule\vbox to 7pt
{\hrule width 6pt\vfill\hrule}\vrule}}}\vskip 0.5cm}

\def\qm{{\widetilde{M}}}
\def\SL{{\bf SL}}
\def\GL{{\bf GL}}
\def\el{{l}}

\let\le=\leq
\let\ge=\geq

\date{\today}
\title{Hyperdifferential properties of \\ Drinfeld quasi-modular forms}
\author{V. Bosser and F. Pellarin}

\begin{document}

\maketitle

{\bf Abstract.} {\footnotesize This article is divided in two parts.
In the first part we endow a certain ring of ``Drinfeld quasi-modular forms''
for $\GL_2(\FF_q[T])$ (where $q$ is a power of a prime) with a system
of ``divided derivatives" (or hyperderivations). This ring contains Drinfeld
modular forms %for $\GL_2(\FF_q[T])$
as defined by Gekeler in \cite{Ge},
and the hyperdifferential ring obtained should be considered
as a close analogue in positive characteristic of famous Ramanujan's
differential system relating to the first derivatives of the classical
Eisenstein series of weights $2$, $4$ and $6$.
%and the ``false'' Eisenstein series of weight $2$.
In the second part of this article
we prove that, when $q\not=2,3$, if ${\cal P}$ is a non-zero
hyperdifferential prime ideal, then it contains the Poincar\'e series
$h=P_{q+1,1}$ of \cite{Ge}. This last result is the analogue of
a crucial property proved by Nesterenko \cite{Nes} in characteristic zero
in order to establish a multiplicity estimate.}

%\tableofcontents

\section{Introduction}\label{intro}

In \cite{Ge}, Gekeler introduced modular forms for $\GL_2(\FF_q[T])$,
analogues in positive characteristic to the classical modular forms for $\SL_2 
(\ZZ)$; these
functions are frequently called ``Drinfeld modular forms''.
Both Drinfeld and classical modular forms have interesting differential
properties but there is a slight imperfection, as
the rings generated by them are not stable under derivation.

To overcome this problem in the classical case, one idea is to
weaken a bit the definition of modular forms, allowing slightly
more general functional equations, thus leading to the notion of {\em
quasi-modular forms} as in \cite[p.~166]{KZ} and
\cite[Definition 113]{MR}
(\footnote{More implicitly, Quasi-modular forms already appear in several
previous works
by Ramanujan and Rankin.}).

The same can be done in the positive characteristic case, as we will  
see in this article: {\em
Drinfeld quasi-modular forms} can be defined, and the ring generated by
these functions both contains Drinfeld modular forms and is stable
under derivation.

To what extent do classical quasi-modular  
forms
and Drinfeld quasi-modular forms resemble  
each other?
This question is the basis for our article.
We will find that, apart from some
differences arising from the unequal characteristics,
the classical theory and Drinfeld's theory of quasi-modular forms bear
striking similarities.

\medskip

\noindent {\bf The classical framework.}
Before we describe our results, it is convenient to present the
aspects of the classical theory whose analogues in positive characteristic
will be discussed in this text.

A holomorphic function $f:{\cal H}\rightarrow\CC$ (where ${\cal H}$
is the complex upper half-plane) is called {\em quasi-modular form}
of weight $w\in\NN$ and depth at most $l\in\NN$ for
$\Gamma_\QQ:=\SL_2(\ZZ)$ if there exist functions $f_0,\ldots,f_l$,
which are holomorphic
both on ${\cal H}$ and ``at infinity'', and which satisfy the following
condition: for all $\gamma=
\sqm{a}{b}{c}{d}\in\Gamma_\QQ$ and all $z\in{\cal H}$,
\begin{eqnarray*}
f\left(\frac{az+b}{cz+d}\right)&=&(cz+d)^w\sum_{r=0}^lf_r(z)\left
(\frac{c}{cz+d}\right)^r.\\
\end{eqnarray*}
After Kaneko and Zagier
(cf. \cite[Proposition~1]{KZ}, see also \cite[Proposition~124]{MR}),
the ring generated by quasi-modular forms for $\Gamma_\QQ$ is also
the ring of functions $\CC[E_2,E_4,E_6]$,
where $E_2,E_4,E_6$ are the classical Eisenstein series of weights
$2,4$ and $6$ respectively.

Recall that $E_4$ and $E_6$ are modular forms (hence,
quasimodular forms of depth $0$),
but that $E_2$ is a non-modular quasi-modular form of depth
$1$. Since the derivative of a quasi-modular
form of weight $w$ and depth $\le l$ is again a quasi-modular form, of
weight $w+2$ and depth $\le l+1$,
and since the functions $E_2,E_4,E_6$ are algebraically
independent over $\CC(q)$ by a classical result due to Rankin and Mahler
(see the contribution of Bertrand in \cite{NP}), $\CC[E_2,E_4,E_6]$
is thus isomorphic to a polynomial ring in three indeterminates, endowed
with a structure of a differential,
filtered (by the depth), graded ring.

Explicitly, if we set $D=(2\pi i)^{-1}d/dz$,
the derivation on this ring is described by
the following formulas (due to Ramanujan):

\begin{equation}\label{eq:systeme_differentiel0}
\left\{\begin{matrix}
DE_2  & = & \frac{1}{12}(E_2^2-E_4 )\\
DE_4  & = & \frac{1}{3}(E_2E_4 -E_6 )\\
DE_6  & = & \frac{1}{2}(E_2E_6 -E_4 ^2).
\end{matrix}\right.
\end{equation}

Among several basic questions that we may ask about differential rings,
one is to determine their
$D$-differential prime ideals (recall that an ideal ${\cal I}$ is said
{\em $D$-differential} if $D {\cal I}\subset {\cal I}$).
In \cite{PF}, the $D$-differential prime ideals of $\CC[E_2,E_4,E_6]$
were all determined
applying elementary commutative algebra. 
There are many others, which we do not describe here,
but they all contain $\Delta$.
This latter property (with a slight technical hypothesis that we
skip here) was used by Nesterenko \cite{Nes} to
prove a {\em multiplicity estimate}, which was a crucial tool
in the proof of his famous theorem on the algebraic independence of
values of Eisenstein series: if $z\in{\cal H}$,
at least three of the four numbers
$e^{2\pi{\rm i}z},E_2(z),E_4(z),E_6(z)$ are algebraically independent over
$\QQ$ (see \cite{Nes}, \cite[Chapter 3]{NP}, see also \cite[Th\'eor\`eme 1.1]{BV}).

\medskip

\noindent {\bf Drinfeldian framework.} After this survey
on quasi-modular forms for the
group $\Gamma_\QQ$,
we are ready to introduce the main results of this article.

Let $T$ be an indeterminate, let $A=\FF_q[T]$ be the polynomial ring
over
a finite field $\FF_q$, and let $K=\FF_q(T)$ be its field of fractions.
Let $K_\infty=\FF_q((1/T))$ be the
completion of $K$ at the infinite prime, and let $C$ be the
completion of an algebraic closure
of $K_\infty$. We will denote by $\Omega$
the set $C\setminus K_\infty$ and, as usual, by
$p$ the characteristic of $C$.

It is well known that $\Omega$ is a connected admissible open
subspace of the rigid analytic space  $\PP_1(C)$. The group
$\Gamma_K:=\GL_2(A)$ acts discontinuously on $\Omega$ by homographies
\[\gamma=\sqm{a}{b}{c}{d}\in\Gamma_K,\quad z\in\Omega,\quad \gamma(z)=
\frac{az+b}{cz+d}.\]
The quotient space can also be given a
rigid structure, and it is natural to consider ``modular
forms'' with the straightforward definition.

The foundations of the theory  of modular forms on $\Omega$ that we
need here are essentially contained in \cite{Ge}. Recall that
since $\Gamma_K$ has the non-trivial character $\det:\Gamma_K
\rightarrow\FF_q^*$,
modular forms naturally have a {\em type} (\footnote{It is a class in
$\ZZ/(q-1)\ZZ$.}), in addition to a weight.
The subring of the ring of holomorphic functions on $\Omega$
generated by modular forms is $C[g,h]$, where $g$ and $h$
are two special modular forms, which are defined in \cite{Ge}.
Here, 
$g$ is the Eisenstein series of weight $q-1$ and type $0$
normalized in a suitable way
(\footnote{In \cite{Ge}, the notation $g$ changes in the
course of the text. In our text we use the notation $g$ for the function
$g_{\text{new}}$ of \cite{Ge} (p.~683).}),
and $h$ is the Poincar\'e series $P_{q+1,1}$ of weight $q+1$ and
type $1$ (example (5.11) of \cite{Ge}).

Our first task (Sections \ref{paragraphe2} to \ref{calculs}) will be
to introduce {\em Drinfeld quasi-modular forms}
for the group $\Gamma_K$ and to study how they behave under
the action of the so-called ``divided derivatives''.

Drinfeld quasi-modular forms will have weights, types and depths, the zero
depth corresponding to the case of modular forms. A useful result
we will then prove is a structure theorem for these quasi-modular forms,
similar to the one of Kaneko and Zagier quoted above. Let
$E$ be the ``false'' Eisenstein series as defined in \cite[p. 686]{Ge}.
With our definition, this function will be a quasi-modular form of
weight $2$, depth $1$ and type $1$.
Denote by $\widetilde{M}_{w,m}^{\leq l}$
the $C$-vector space of quasi-modular forms of weight $w$, type $m$
and depth $\leq l$, by $\widetilde{M}$ the ring generated by
all quasi-modular forms (see \S~\ref{paragraphe2} for the precise
definitions), and by $C[E,g,h]_{w,m}^{\le\el}$
the subspace of $C[E,g,h]$ generated by the monomials
$E^{\alpha}g^{\beta}h^{\gamma}$ satisfying the conditions
$2\alpha+\beta(q-1)+\gamma(q+1)=w$,
$\alpha+\gamma\equiv m\pmod{q-1}$ and $\alpha\le\el$.

In Section~\ref{paragraphe2} we prove:

\begin{Theo_intro}\label{theo1}
The functions $E$, $g$, $h$ are algebraically independent over $C$,
and we have
$$ \qm_{w,m}^{\le\el}=C[E,g,h]_{w,m}^{\le\el}\qquad
{\text and}\qquad \qm=C[E,g,h].$$
\end{Theo_intro}

Before introducing the differential results of this paper, let us recall that
Gekeler proved that the three functions $E$, $g$ and $h$
satisfy a system of differential equations, namely (\footnote{See
\cite{Ge}: the last equation corresponds to the definition of $E$
(p. 686), while the first relation corresponds to (8.6)
and the second relation is equivalent to Theorem 9.1.}):
\begin{equation}\label{system}
\left\{\begin{array}{ccl}
D_1E& = & E^2\\
D_1g& = & -(Eg+h)\\
D_1h& = & Eh,
\end{array}\right.
\end{equation}
where $D_1$ denotes the operator $(-\overline{\pi})^{-1}d/dz$,
$\overline{\pi}$ being a fixed fundamental period of the Carlitz
module.
Thus, $D_1$ defines a derivation on the ring $C[E,g,h]$.
At first glance, the system (\ref{system}) seems to be an analogue of
(\ref{eq:systeme_differentiel0}). But it turns out that the kernel of
the derivation $D_1$ is far too ``big'' (\footnote{It can be proved
that it is equal to the ring
\[C[E^p,g^p,h^p,Eg+h,Eh^{p-1},E^2h^{p-2},\ldots,E^{p-1}h].\]}), and
it is hopeless to get results about differential
ideals analogous to those mentioned above in the classical case.
For instance, the kernel of $D_1$ contains infinitely
many irreducible elements, and so the differential ring $C[E,g,h]$
has infinitely many $D_1$-differential principal prime ideals,
unlike the characteristic zero case.

The first important objective of this article is to show
that in fact, the ring $C[E,g,h]$
is not only stable under the derivation $D_1$, but also under
all the ``divided derivations''.

Let $f:\Omega\rightarrow C$ be a holomorphic function and
$z\in\Omega$ be a point. We define the \emph{divided derivatives} or
\emph{hyperderivatives}
$({\cal D}_nf)(z)$ of $f$ at $z$ ($n\ge 0$) by the formula
$$
f(z+\varepsilon)=\sum_{n\ge 0}({\cal D}_nf)(z)\varepsilon^n,\
\varepsilon\in C,\ \ |\varepsilon|\ \text{small.}
$$
This defines holomorphic functions ${\cal D}_nf:\Omega\rightarrow C$.
In order to preserve $K$-rationality properties, define further
$D_n:=1/(-\overline{\pi})^n{\cal D}_n$ ($n\ge 0$).
We then have the following result.

\begin{Theo_intro}\label{theo2}
Let $n\ge 0$ be an integer. The operators $D_n$ induce $C$-linear operators
$$ D_n:C[E,g,h]\rightarrow C[E,g,h] $$
such that, for all $(w,l,m)\in\NN\times\NN\times\ZZ/(q-1)\ZZ$,
\begin{equation}
D_n\bigr(K[E,g,h]_{w,m}^{\le\el}\bigl)\subset
K[E,g,h]_{w+2n,m+n}^{\le\el+n}.\label{inclusion}
\end{equation}
\end{Theo_intro}

This Theorem will be proved in Section~\ref{section_hyp}
as a consequence of Proposition \ref{propo}.
It is the {\em hyperdifferential} ring $C[E,g,h]$ (that is, the ring equipped
with the operators $D_1,D_2,\ldots$) which should be interpreted as
the analogue in the Drinfeldian framework of the differential ring
$\CC[E_2,E_4,E_6]$.

By Theorem~\ref{theo2}, any hyperderivative of a polynomial
in the functions $E$, $g$, $h$ is again a polynomial in these
functions. However, the operators $D_n$ behave erratically and the problem
of finding explicit formulas for $D_nE$, $D_ng$, and $D_nh$ when $n\ge p$
is not a trivial one, in contrast with the classical case, where it suffices
to iterate the operator $D$ (see for example the formula (25) of \cite{SD}).
Some general results concerning
the hyperderivatives $D_nE$, $D_ng$ and $D_nh$ can
nevertheless be obtained when $n$ has a particular form. For example,
in Proposition~\ref{derivees_mod_modulaires} we determine
$D_{p^k}E$ and $D_{p^k}g$ for all $k$, up to a modular form
(see Section~\ref{uptomodular} for other results in the same vein).
Also, we give explicit formulas in Section~\ref{calculs}
for $D_{p^k}E$, $D_{p^k}g$ and $D_{p^k}h$ when $p^k\le q^2$
(Theorem~\ref{theo_calculs}).
As we will see, these formulas allow, at least in principle, to compute
$D_nf$ for any $f\in C[E,g,h]$ and any $n\le q^2$.
For instance, in the particular case $q=p$, these formulas (which
then determine $D_nf$ even for all $n\le p^3-1$) read

\begin{equation}\label{systemep}
\left\{\begin{array}{ccl}
D_p E&=&E^{p+1}+\frac{1}{d_1}h^2\\
D_p g&=&E^pg\\
D_p h&=&2E^ph-\frac{1}{d_1}gh^2,
\end{array}\right.
\end{equation}
\begin{equation}\label{systemep2}
\left\{\begin{array}{ccl}
D_{p^2}E&=&E^{p^2+1}+\frac{1}{d_1^p}g^{p-1}h^{p+1}+\frac{1}{d_2}g^{2p}
h^2\\
D_{p^2}g&=&E^{p^2}g-\frac{d_1}{d_2}g^{p+1}h^p+
\frac{d_2-d_1^{p+1}}{d_1^{p-1}d_2}h^{2p-1}\\
D_{p^2}h&=&E^{p^2}h+\frac{1}{d_1^{p-1}}E^pg^{p-1}h^p-\frac{1}{d_2}
g^{2p+1}h^2-\frac{d_1^{p+1}+d_2}{d_1^pd_2}g^ph^{p+1},
\end{array}\right.
\end{equation}
where $d_1=[1]$ and $d_2=[2][1]^q$, the notation
$[i]$ being defined as usual by $[i]=T^{q^i}-T$.

\medskip

The second aspect of the theory that we want to develop is, in view of
Neste\-renko's results quoted before, the study of the ideals which are
stable by all the operators $D_n$.
If ${\cal I}$ is an ideal of the ring $C[E,g,h]$, we say that
${\cal I}$ is \emph{hyperdifferential},
if it is stable under all the operators $D_n$,
that is, if $D_n({\cal I})\subset {\cal I}$ for all $n\ge 0$. 
Hyperdifferential ideals in $C[E,g,h]$ are natural analogues of
$D$-differential ideals of $\CC[E_2,E_4,E_6]$.

Non-trivial hyperdifferential ideals
of $C[E,g,h]$ exist: for example, in Proposition
\ref{hstable} we will see that the principal ideal $(h)$ is hyperdifferential,
in analogy with the ideal $(\Delta)$ in the classical case.
The next important result we will prove is that in fact, when $q\not=2,3$,
this ideal is the \emph{only} non-zero hyperdifferential principal
prime ideal (again as in the classical case).

More precisely, let us write $F=C(E,g,h)$. The hyperderivatives
$D_n$ ($n=1,2,\ldots$) extend to $F$ in a unique way, we denote
these extensions by $D_n$ again. We have:

\begin{Theo_intro}
Assume that $q\not\in\{2,3\}$.
Let $f\in F^\times$ be such that for all $r\geq 0$
$(D_{r}f)/f\in C[E,g,h]$. Then, there exists $n\in \ZZ$ and
$c\in C^\times$ such that $f=ch^n$.
\label{theoreme_h}\end{Theo_intro}

Note that our proof of Theorem~\ref{theoreme_h} really requires
that $q\not=2,3$, and that we do not know if the result remains true
when $q\in\{2,3\}$.
This Theorem is an analogue of Lemma 5.2 of \cite[Chapter~10]{NP}
(or Lemma 4.26 of \cite{BV}), but its proof is different, and
much more difficult. It requires two steps, contained in Section
\ref{section:groups}. The first step consists
in showing that if $f$ is as in
Theorem~\ref{theoreme_h}, there exists a $p$-adic integer
$\sigma=s_1+s_2p+\cdots$ such that for all $k\geq 0$,
$fh^{-(s_1+s_2p+\cdots+s_{k+1}p^k)}\in\ker D_1\cap\cdots\cap\ker D_{p^k}$.
Then, with the help of the expansions at infinity of the functions,
we show that $\sigma\in\ZZ$,
which yields the result since $\cap_{k\ge 0}\ker D_{p^k}=C$.

Somewhat, this approach already appears in the paper
\cite{MP} by Matzat and van der Put. But these authors work with
{\em linear} hyperdifferential modules over a field, while in our context,
the theory is essentially non-linear. This means that our proofs, even though
initially inspired by  the work \cite{MP}, ultimately depend very
slightly on it. We should also add that in the proof of
Theorem~\ref{theoreme_h}, it is not enough to know the formulas
(\ref{system}) (or even the formulas of Section~\ref{calculs}).
Thus, in this part of our work, the
formalism of Kaneko and Zagier becomes very helpful.

As in characteristic zero, one can use Theorem~\ref{theoreme_h} to prove a
result for arbitrary non-zero hyperdifferential prime ideals. One gets:

\begin{Corollaire_intro}
Assume that $q\not\in\{2,3\}$.
If ${\cal P}$ is a non-zero hyperdifferential prime ideal of $C[E,g,h]$,
then $h\in{\cal P}$.
\label{theoreme_hh}\end{Corollaire_intro}

Corollary~\ref{theoreme_hh} is proved in
Section~\ref{section:non_principal}. It is
an analogue of Proposition 5.1 of
\cite[Chapter~10]{NP} (see also Theorem 4.25 of \cite{BV}).
To prove it, we see from Theorem~\ref{theoreme_h} that
we only have to treat the case when
the height of the prime ideal ${\cal P}$ is $\geq 2$. In characteristic zero,
to handle this case the original argument by Nesterenko was
a generalisation of an idea of Siegel for the classification of the
algebraic solutions of Bessel and Riccati differential equations.
This method does not seem to work in our case; to avoid it, we apply a
trick appealing to Rankin brackets as in \cite{PF}, which also
works in positive characteristic.

Finally, in Section~\ref{section:examples} (Theorem~\ref{stable_ideals}),
we apply Theorem~\ref{theoreme_h} and Corollary~\ref{theoreme_hh} to obtain
a full description of all the hyperdifferential ideals of $C[E,g,h]$,
analogous to Theorem 1.2 of \cite{PF}.

\medskip

\noindent \emph{Remarks about Nesterenko's Theorem.} Since it was
one of the first motivations of our work,
we cannot finish this introduction without a discussion around a possible
analogue in the Drinfeldian framework of the above quoted Nesterenko's
Theorem. Let $t(z)$ be the usual ``parameter at infinity''
(see Section~\ref{paragraphe2} for the definition). It seems natural to
state the following Conjecture (compare with \cite[Conjecture 1]{Dion}).

\medskip

\noindent\textbf{Conjecture.} {\em For all
$z\in\Omega$, at least three of the four numbers $$t(z),E(z),g(z),h(z)$$
are algebraically independent over $K$.}

\medskip

In the direction of this Conjecture, an important progress would be to
prove an analogue of Nesterenko's multiplicity estimate (see
\cite[Theorem 3]{Nes}, see also \cite[Th\'eor\`eme 2.9]{BV}). Here too, we
could make a Conjecture but we skip it as it merely consists of substituting
$q$ and the $q$-expansions of $E_2,E_4,E_6$ with $t$ and the $t$-expansions
of $E,g,h$ in Nesterenko's multiplicity estimate (see
Section~\ref{paragraphe2} for the precise definition of $t$-expansion).

As already mentioned, the analogues of Theorem~\ref{theoreme_h} and
Corollary~\ref{theoreme_hh} in characteristic zero play a decisive role
in the proof of Nesterenko's estimate (see \cite[Remarque 4.5 and Th\'eor\`eme 4.22]{BV}).
However,
although the results we obtain here are very similar to those in characteristic
zero, they do not suffice to reach a multiplicity estimate: just as
the proof of Theorem~\ref{theoreme_h} is much more difficult and need
new arguments, it seems that a proof of the multiplicity estimate
in characteristic $p$ will also require further ideas (we omit the details about the several occurring obstructions).

Consequently, the above Conjecture seems to be still
out of reach. It is worth noting that such a situation is
quite atypical: indeed, transcendence (and algebraic independence) theory in
the Drinfeldian framework has now strongly been developed,
and most of the classical transcendence results have a known analogue
in positive characteristic. Recently, the theory has even gone far beyond
the classical one in some cases (\cite{abp}, \cite{papa}, \cite{changyu}).

Let us notice that, just as for Nesterenko's Theorem, the Conjecture does not
seem to follow from Grothendieck-like period conjectures in positive
characteristic, but only from some variants of Andr\'e's conjecture. Hence,
it is unclear how it could be handled with the methods of \cite{papa}.

The results we obtain in this text hopefully constitute
a first step towards the Drinfeldian analogue of Nesterenko's
multiplicity estimate. We hope taking others (and maybe all?) in a
subsequent paper.

\section{Basic properties of quasi-modular forms.\label{paragraphe2}}

The aim of this Section is to introduce Drinfeld quasi-modular forms
and to prove some of their fundamental properties, in particular
Theorem~\ref{theo1}.

Recall that in the previous Section we have denoted by $\overline{\pi}$
a fixed fundamental period of the Carlitz module. In the whole text, we will
further denote by $e_{{\mathcal C}}:C\rightarrow C$ the Carlitz exponential,
and by $t:\Omega\rightarrow C$ the usual ``parameter at infinity'', that is,
$t(z)=1/e_{{\mathcal C}}(\overline{\pi}z)$. We have $\ker e_{{\mathcal C}}=
\overline{\pi} A$ and
\begin{equation}\label{exponential}
e_{{\mathcal C}}(z)=\sum_{n\ge 0}\frac{z^{q^n}}{d_n},
\end{equation}
where $d_i:=[i][i-1]^q\cdots[1]^{q^{i-1}}$ and $[i]:=T^{q^i}-T$
(see Section~4 of \cite{Ge}).

Following the terminology of \cite{Ge} and \cite{vdp}, we will call
{\em holomorphic function} in $\Omega$ an analytic function on $\Omega$
in the rigid analytic sense.
Holomorphic functions are thus {\em locally analytic} in the sense of \cite{US}
(but the converse is false).

The {\em imaginary part} of an element $z\in C$
(see Section 5 of \cite{Ge}) will be denoted by $|z|_i$. If
$f:\Omega\rightarrow C$ is holomorphic in $\Omega$ and $A$-periodic
({\em i.e.} $f(z+a)=f(z)$ for all $a\in A$ and $z\in \Omega$), we will say
that it is {\em holomorphic at infinity} if
$f(z)$ is the value, for $|z|_i$ sufficiently large, of a converging series
\begin{equation}
f(z)=\sum_{i\geq 0}c_it(z)^i,\label{eq:texp}
\end{equation}
where $c_0,c_1,\ldots$ are elements of $C$.
These coefficients uniquely determine $f$ and we will often write, by abuse of
notation, $f=\sum_{i\geq 0}c_it^i$.
The series (\ref{eq:texp}) is the {\em $t$-expansion} of $f$.

In the whole text, we will denote by ${\cal O}$ the ring of holomorphic
functions $f:\Omega\rightarrow C$ which are $A$-periodic and
holomorphic at infinity. This ring contains in particular the modular
forms for $\Gamma_K$ (by Definition 5.9 of \cite{Ge}).
If $f$ is a non-zero element of
${\cal O}$, we will denote by $\nu_\infty(f)$ its {\em order at infinity}:
this is the smallest integer $i$ such that $c_i\not=0$ in the
$t$-expansion (\ref{eq:texp}). The map $\nu_{\infty}$ extends
to a map $\nu_\infty:\mbox{Frac}({\cal O})\setminus\{0\}\rightarrow\ZZ$ in the obvious way.

\begin{Definition}
Let $w\ge 0$ be an integer and $m\in\ZZ/(q-1)\ZZ$.
A holomorphic function $f:\Omega\rightarrow C$ is called
\emph{quasi-modular form of weight $w$ and type $m$} if there exist
functions $f_0,\ldots,f_{\el}\in{\cal O}$ such that, for all $z\in 
\Omega$
and $\gamma=\begin{pmatrix} a&b\\ c&d\end{pmatrix}\in\Gamma_K$:
\begin{equation}
f(\gamma(z))=(cz+d)^w(\det\gamma)^{-m} \sum_{i=0}^{\el} f_i(z)
\bigl(\frac{c}{cz+d}\bigr)^i.\label{QMcondition}
\end{equation}
\end{Definition}

The type $m$ of a quasi-modular form is a class in $\ZZ/(q-1)\ZZ$.
However, to avoid heavy notations, we will sometimes
use the same notation $m$ for any
representant of the class, thus considering $m$ as an integer.
This abuse of notation will not lead to confusion.

\medskip

\noindent \textbf{Remarks.} (i). Let $f$ be a quasi-modular form of weight
$w$ and type $m$, and
let $f_0,\ldots,f_{\el}\in{\cal O}$ be as in the definition.
It is easy to see that if $f\not=0$, then the weight $w$, the type $m$
and the polynomial $\sum_{i=0}^{\el}f_iX^i\in{\cal O}[X]$
are uniquely determined by $f$. Indeed, if the function $f$ were also a  
quasi-modular form of weight $w'\ge w$
and type $m'$, and if an identity analogous to (\ref{QMcondition}) held
with functions $g_0,\ldots,g_{\el'}\in {\cal O}$ at the place of $f_0, 
\ldots,f_\el$, then we would have
\begin{equation}
c^{w-w'}(\det\gamma)^{m'-m}A_z(\frac{c}{cz+d})-B_z(\frac{c}{cz+d})=0
\label{bbb1}
\end{equation}
for all $z\in\Omega$ and all $\gamma\in\Gamma_K$ with $c\not=0$,
where
$$
A_z=X^{w'-w}\sum_{i=0}^{\el}f_i(z)X^i\in C[X]\quad and\quad
B_z=\sum_{i=0}^{\el'}g_i(z)X^i\in C[X].
$$
Choosing now $\gamma=\begin{pmatrix} 1&d-1\\ 1&d\end{pmatrix}\in 
\Gamma_K$
in (\ref{bbb1}) (with $d\in A$ arbitrary), we obtain
$A_z=B_z$ for all $z\in\Omega$. By fixing an element $z\in\Omega$
such that $A_z\not=0$, choosing $c\in A$ with $\deg c\ge 1$
such that $A_z(c/(cz+d))\not=0$, and taking
$\gamma=\begin{pmatrix} \lambda & 0 \\ c & 1 \end{pmatrix}$
with $\lambda\in\FF_q^*$ in the equation (\ref{bbb1}), we get
$$ \lambda^{m'-m}=c^{w'-w}\quad \text{for all\ }\lambda\in\FF_q^*.$$
This implies $w=w'$ and $m = m'$ in $\ZZ/(q-1)\ZZ$. The equality $A_z=B_z$
for all $z$ then yields the equality between the two polynomials
$\sum_{i=0}^{\el}f_iX^i$ and  $\sum_{i=0}^{\el'}g_iX^i$, as claimed.

In virtue of this remark, if $f$ is any non-zero quasi-modular
form, the polynomial
\begin{equation*}
P_f:=\sum_{i=0}^{\el}f_iX^i
\end{equation*}
of ${\cal O}[X]$ is well defined. We call it the
{\em associated polynomial} of $f$.
The degree of $P_f$ is called the \emph{depth} of $f$ and is usually
denoted by $\el$; another notation is $l(f)$. If $f=0$, we  
agree that it is a quasi-modular form
of weight $w$, type $m$ and depth $\el$ for all $w,m,\el$, and we set
$P_f:=0$.

We denote by $\qm_{w,m}^{\el}$ (\resp\ by $\qm_{w,m}^{\le\el}$)
the set (\resp\ the $C$-vector space) of quasi-modular forms of weight
$w$, type $m$ and depth $\el$ (\resp\ depth $\le\el$). Note that
$\qm_{w,m}^{\el}$ is not a vector space for $\el\ge 1$.
We will further denote by $\qm_{w,m}$ the $C$-vector space
of quasi-modular forms of weight $w$ and type $m$, and by
$\widetilde{M}$ the $C$-algebra $\widetilde{M}=\sum_{w,m}\qm_{w,m}$.
In fact, this sum is a direct sum (see Proposition~\ref{prop1}
below).

\medskip

\noindent (ii). If $f$ is an element of $\qm_{w,m}^{\el}$, then the
constant term $f_0$ of the associated polynomial $P_f$ is
necessarily equal to $f$, as
follows from the formula (\ref{QMcondition}) by choosing for $\gamma$
the identity  matrix.

\medskip

\noindent (iii).
It is clear from the definition that a modular form $f$ of weight $w$  
and
type $m$ is an element of $\qm_{w,m}^0$ with associated polynomial
$P_f=f$. Conversely, Remark (ii) shows
that any element $f$ of $\qm_{w,m}^0$ is a modular form of weight
$w$ and type $m$, thus we have:
\begin{equation}
M_{w,m} = \qm_{w,m}^0,\label{bbb0}
\end{equation}
where $M_{w,m}$ denotes the $C$-vector space of
Drinfeld modular forms of
weight $w$ and type $m$ (equal to $M_w^m$, following the notation
of \cite{Ge}).

\medskip

\noindent (iv). If $f_1\in\qm_{w_1,m_1}^{\el_1}$ and
$f_2\in\qm_{w_2,m_2}^{\el_2}$,
it is obvious from the definitions
that $f_1f_2\in\qm_{w_1+w_2,m_1+m_2}^{\el_1+\el_2}$ and that
$P_{f_1f_2}=P_{f_1}P_{f_2}$. Moreover, $P_{\lambda f_1}=\lambda P_{f_1}$
for $\lambda\in C$, and if $w_1=w_2$, $m_1= m_2$, then
$f_1+f_2\in \qm_{w_1,m_1}^{\le\max\{\el_1,\el_2\}}$
with $P_{f_1+f_2}=P_{f_1}+P_{f_2}$.

\medskip

An important example of a quasi-modular form which is not a modular  
form is
furnished by the ``false'' Eisenstein series $E(z)$,
which satisfies (see \cite[formula (8.4)]{Ge})
\begin{equation}
E(\gamma(z))=(cz+d)^2(\det\gamma)^{-1}\left(E(z)-\frac{c}{\overline 
{\pi}(cz+d)}\right)
\label{bbb3}
\end{equation}
for all $z\in\Omega$ and all $\gamma\in\Gamma_K$.
In virtue of the results of Section 8 of \cite{Ge}, we have $E\in {\cal O}$.
Thus, $E\in\qm_{2,1}^1$ and $P_E=E-\overline{\pi}^{-1}X$.

\begin{Proposition}\label{prop1}
We have:
$$ \widetilde{M}=
\bigoplus_{\genfrac{}{}{0pt}{1}{w\ge 0}{m\in\ZZ/(q-1)\ZZ}}
\qm_{w,m}.$$
\end{Proposition}

\noindent\emph{Proof.} The result can be proved exactly as in
\cite[proof of Lemma~16]{MR} and is thus left to the reader. \CVD

Throughout this text,
we will denote by $C[E,g,h]_{w,m}^{\le\el}$
the space of elements of $C[E,g,h]$ of weight $w$, type $m$ and depth
$\le\el$,
that is, the $C$-vector space of polynomials in $E$, $g$, $h$ which  
are linear
combinations with coefficients in $C$ of monomials
$E^{\alpha}g^{\beta}h^{\gamma}$
with $2\alpha+\beta(q-1)+\gamma(q+1)=w$, $\alpha+\gamma\equiv m\pmod{q-1}$
and $\alpha\le\el$.

\medskip

To prove Theorem 1, we will follow a method similar to the one used in
\cite{MR} in the complex case. We need four Lemmata.
In the next Lemma, we will say that a polynomial $P\in C[Y_1,Y_2,Y_3]$
is {\em isobaric} of weight $w\in\NN$ if it is the
sum of monomials of the form $\lambda Y_1^{\alpha}Y_2^{\beta}Y_3^ 
{\gamma}$,
where $\lambda\in C$ and $\alpha$, $\beta$, $\gamma$ satisfy
$2\alpha+(q-1)\beta+(q+1)\gamma=w$. Thus, an isobaric element of
$C[Y_1,Y_2,Y_3]$ is nothing but a homogeneous element
for the graduation defined associating to $Y_1,Y_2,Y_3$ the weights  
$2,q-1,q+1$.

If all the monomials above are moreover such that $\alpha+\gamma
\equiv m\pmod{q-1}$, then we
will say that $P$ is of type $m\in\ZZ/(q-1)\ZZ$.

\begin{Lemme}\label{lemme1}
Let $P\in C[Y_1,Y_2,Y_3]$ be a non-zero polynomial isobaric of weight  
$w$ and type $m$.
Then the function $P(E,g,h)$ is a quasi-modular form of weight $w$,
type $m$ and depth $\el=\deg_{Y_1}P$.
\end{Lemme}

\noindent\emph{Proof.} Write
$$ P = \sum_{\alpha=0}^{\el}p_{\alpha}(Y_2,Y_3)Y_1^{\alpha},$$
where $p_{\el}\not=0$ and $p_{\alpha}$ is an isobaric polynomial
of weight $w-2\alpha$ and type $m-\alpha$ ($0\le\alpha\le\el$).
Each function $p_{\alpha}(g,h)E^{\alpha}$ is clearly an element of
$\qm_{w,m}^{\alpha}$. But it is plain that the sum of two non-zero  
elements
of $\qm_{w,m}$ of different depths $\alpha_1$, $\alpha_2$
is again an element of $\qm_{w,m}$ of depth $\max\{\alpha_1,\alpha_2\}$.
It follows at once that $P(E,g,h)\in\qm_{w,m}^{\deg_{Y_1}P}$. \CVD

\begin{Lemme}\label{lemme2}
The functions $E$, $g$ and $h$ are algebraically independent over $C$.
\end{Lemme}

\noindent\emph{Proof.} Let $P\in C[Y_1,Y_2,Y_3]$
be such that $P(E,g,h)=0$. We can write
$$
P=\sum_{w=0}^d\ \sum_{m\in\ZZ/(q-1)\ZZ}P_{w,m},$$
where $P_{w,m}$ is an isobaric polynomial of weight $w$ and type $m$.
From $P(E,g,h)=0$ and from Proposition~\ref{prop1} it follows
that $P_{w,m}(E,g,h)=0$ for each pair $(w,m)$. We fix $w$ and $m$
and write $$P_{w,m} =\sum_{\alpha\ge 0}p_{\alpha}(Y_2,Y_3)Y_1^{\alpha}, $$
where $p_{\alpha}\in C[Y_2,Y_3]$ is isobaric of weight $w-2\alpha$ and type
$m-\alpha$.
Suppose that there is an index $\alpha_0\ge 1$ such that
$p_{\alpha_0}(g,h)\not=0$, and take $\alpha_0$ maximal. We have
\begin{equation}
p_{\alpha_0}(g,h)E^{\alpha_0}=
-\sum_{0\le\alpha<\alpha_0}p_{\alpha}(g,h)E^{\alpha}.
\label{bbb2}
\end{equation}
Now, $p_{\alpha_0}(g,h)E^{\alpha_0}$
is a non-zero quasi-modular form of depth exactly $\alpha_0$, and
the right-hand side of (\ref{bbb2}) is a quasi-modular form of depth  
at most
$\alpha_0-1$ (Lemma~\ref{lemme1}). By uniqueness of the depth,
this is impossible.
Hence we deduce that $p_{\alpha}(g,h)=0$ for all  $\alpha$.
Since now $g$ and $h$ are known to be algebraically independent over
$C$, we get $p_{\alpha}=0$ for all $\alpha$, and thus $P_{w,m}=0$.
This being true for all $(w,m)$, we have $P=0$, which completes the
proof of the Lemma. \CVD

\begin{Lemme}\label{lemme3}
Let $f\in\qm_{w,m}^{\el}$, and let $P_f=\sum_{i=0}^{\el}f_iX^i$
be its associated polynomial. Then, for all $i$, $0\le i\le\el$,
the function $f_i$ is an element of the space
$\qm_{w-2i,m-i}^{\le\el-i}$, whose associated
polynomial $P_{f_i}$ is explicitly given by
$$
P_{f_i}=\sum_{j=i}^{\el}\binom{j}{i} f_jX^{j-i}.
$$
In particular, we have $f_{\el}\in M_{w-2\el,m-\el}$.
\end{Lemme}

\noindent\emph{Proof.} The same proof as in \cite[Lemma 119]{MR} works.
First, an easy computation shows that for
$N=\begin{pmatrix} \alpha&\beta\\ \gamma&\delta\end{pmatrix}
\in\Gamma_K$, $z\in\Omega$ and $i\in\NN$ with $0\le i\le\el$,
the following equality holds:
$$
f_i(z)=\sum_{j=i}^{\el}\binom{j}{i}(-\gamma)^{j-i}
(\gamma(z)+\delta)^{i+j-w}(\det N)^{m-j}f_j(Nz).
$$
Substituting $z$ by $L(z)$, where $L=\begin{pmatrix}a&b\\
c&d\end{pmatrix}:=N^{-1}$, then gives the following transformation
formula for $f_i(L(z))$:
$$
f_i(L(z))=(cz+d)^{w-2i}(\det L)^{-(m-i)}
\sum_{j=i}^{\el}\binom{j}{i} f_j(z)\bigl(\frac{c}{cz+d}\bigr)^{j-i}.
$$
The assertions of the Lemma follow at once. \CVD

\begin{Lemme}\label{lemme4}
For all integers $w\ge 0$, $\el\ge 0$ and for all $m\in\ZZ/(q-1)\ZZ$,
we have:
$$ \qm_{w,m}^{\le\el}=\bigoplus_{0\le i\le\el}M_{w-2i,m-i}E^i.$$
\end{Lemme}

\noindent\emph{Proof.} We already know from Lemma~\ref{lemme2} that
the sum
$$\sum_{0\le i\le\el}M_{w-2i,m-i}E^i$$
is a direct sum, and that it is a subspace of $\qm_{w,m}^{\el}$
(Lemma~\ref{lemme1}). So it remains to prove the opposite inclusion.
We argue by induction on $\el$. If $\el=0$, the result is true (it
is equality (\ref{bbb0})). Suppose now that the result is true for
$\el-1$, where $\el\ge 1$. Let $f\in\qm_{w,m}^{\el}$, $f\not=0$,
and let $P_f=\sum_{i=0}^{\el}f_iX^i$ be its associated polynomial.
Since $P_E=-\overline{\pi}^{-1}X+E$ (see formula (\ref{bbb3})), we have
$P_{E^{\el}}=(-\overline{\pi}^{-1}X+E)^{\el}$. Thus, since
we have $f_{\el}\in M_{w-2\el,m-\el}$ by Lemma~\ref{lemme3},
we get
$$
P_{f_{\el}E^{\el}}=(-\overline{\pi})^{-\el}f_{\el}X^{\el}+
\text{terms of degree} < \el.
$$
Define now $F:=f-(-\overline{\pi})^{\el}f_{\el}E^{\el}$. The  
function $F$
is a quasi-modular form of weight $w$ and type $m$ with associated
polynomial $P_F=P_f-(-\overline{\pi})^{\el}P_{f_{\el}E^{\el}}$,
whose degree is by construction at most $\el-1$. Hence
$F\in\qm_{w,m}^{\le\el-1}$. By induction we can write now
$F=\sum_{i=0}^{\el-1}F_iE^i$, where $F_i\in M_{w-2i,m-i}$ for all $i$,
from which it follows
$$f=F+(-\overline{\pi})^{\el}f_{\el}E^{\el}
\in\bigoplus_{0\le i\le\el}M_{w-2i,m-i}E^i.$$ \CVD

\noindent\emph{Proof of Theorem~\ref{theo1}.}
Since $C[E,g,h]_{w,m}^{\le 0}=M_{w,m}$ (see \cite{Ge}), we obviously have
$$C[E,g,h]_{w,m}^{\le\el} =\bigoplus_{0\le i\le\el}M_{w-2i,m-i}E^i.$$
Therefore Theorem~\ref{theo1} follows from Lemmata \ref{lemme2} and
\ref{lemme4}.
\CVD

Let $L$ be a subfield of $C$ containing $K$ (recall
that $K=\FF_q(T)$). Denote
by $\qm_{w,m}^{\le\el}(L)$ the $L$-vector space of elements of
$\qm_{w,m}^{\le\el}$ having a $t$-expansion at infinity with
coefficients in $L$, and define $L[E,g,h]_{w,m}^{\le\el}$
similarly as before, \ie\ as the $L$-vector space
$C[E,g,h]_{w,m}^{\le\el}\cap L[E,g,h]$. Theorem~\ref{theo1} gives rise
to the following rationality result.

\begin{Corollaire}\label{coroll}
Let $L\subset C$ be any subfield of $C$ containing $K$. Then
$$\qm_{w,m}^{\le\el}(L)=L[E,g,h]_{w,m}^{\le\el}.$$
\end{Corollaire}

\noindent\emph{Proof.} Since the functions $E,g,h$ are elements of
$\qm_{w,m}^{\le\el}(L)$ (the coefficients of their $t$-expansions are
in $A$, see \cite{Ge}), we already have the inclusion
\begin{equation}
L[E,g,h]_{w,m}^{\le\el}\subset\qm_{w,m}^{\le\el}(L).\label{inclus}
\end{equation}
Extending now the scalars from $L$ to $C$, we get
$$C[E,g,h]_{w,m}^{\le\el}\subset\qm_{w,m}^{\le\el}(L)\otimes_L C
\subset \qm_{w,m}^{\le\el}(C)=\qm_{w,m}^{\le\el}.$$
By Theorem~\ref{theo1}, these inclusions between $C$-vector spaces
are in fact equalities, and thus
\begin{align*}
\dim_{L}\qm_{w,m}^{\le\el}(L)=\dim_C\bigl(
\qm_{w,m}^{\le\el}(L)\otimes_L C\bigl)& = \dim_CC[E,g,h]_{w,m}^{\le
\el}\\
& = \dim_LL[E,g,h]_{w,m}^{\le\el}.
\end{align*}
It then follows from (\ref{inclus}) that
$L[E,g,h]_{w,m}^{\le\el}=\qm_{w,m}^{\le\el}(L)$.\CVD

\section{Hyperderivatives and quasi-modular forms\label{section_hyp}}

In this Section, we study the hyperderivatives of quasi-modular forms.
After recalling some basic properties on higher derivations and
hyperderivatives in Section~\ref{basic}, we prove in
Section~\ref{pol_associated} a very important proposition
(Proposition~\ref{propo}),
which gives the connection between the polynomial $P_f$ associated to a
quasi-modular form $f$ and the polynomial $P_{{\cal D}_nf}$ associated to
its $n$-th hyperderivative ${\cal D}_nf$. Using this
result, we can then prove Theorem~\ref{theo2} in Section~\ref{proof_th2},
and the fact that the ideal $(h)$ is hyperdifferential in
Section~\ref{h_hyperstable}. Finally, we determine
in Section~\ref{uptomodular}, up to a modular form, the hyperderivatives
$D_{p^i}E$, $D_{p^i}g$ for all integers $i\ge 0$, and we prove
a similar result for the $D_{p^i}h$'s.

In what follows, we put as usual $\binom{n}{i}=0$ if
$0\le n< i$ and $\binom{n}{i}=(-1)^i\binom{i-n-1}{i}$ if $i\ge 0$, $n<0$.
In particular, we have $\binom{n}{0}=1$ for all $n\in\ZZ$.

We also recall the following fact, which will be often used
in the sequel: if $r$ is any positive integer and
$n=n_0+n_1p^r+\ldots+n_sp^{rs}$, $i=i_0+i_1p^r+\ldots+i_sp^{rs}$ are
two elements of $\NN$ written in base $p^r$, then one has in
$\FF_p$:
\begin{equation}\label{eq:lucas}
\binom{n}{i}=\binom{n_0}{i_0}\binom{n_1}{i_1}\cdots\binom{n_s}{i_s}.
\end{equation}
The formula (\ref{eq:lucas}) was first observed by Lucas in
\cite{Lucas:Binomiaux} (strictly speaking for $r=1$) and
can be easily derived by identifying the coefficient of $X^i$ in both
sides of the following equality of polynomials of $\FF_p[X]$
(we expand of course each factor occurring
in both sides with the help of the binomial formula):
$$(1+X)^n=(1+X)^{n_0}(1+X^{p^r})^{n_1}\cdots(1+X^{p^{rs}})^{n_s}.$$

\subsection{Background on higher derivations and hyperderivatives}\label{basic}

We begin by a short review on iterative higher
derivations, following Matsumura \cite[Section 27]{MH}.

After \cite{MH},
a {\em higher derivation} $(D_n)_{n\in\NN}$ on a $C$-algebra $R$ is a family
of $C$-linear maps from $R$ to itself, with $D_0$ the identity map, such that
if $f,g\in R$, then the following Leibniz's rule holds
\begin{equation}
D_i(fg)=\sum_{r=0}^i(D_rf)(D_{i-r}g).\label{eq:op2}
\end{equation}
In particular, the following formula holds for $f\in R$ and for all
integers $n,k\ge 0$ (see \cite[p.~209]{MH})
\begin{equation}
D_{np^{k}}(f^{p^k})=(D_nf)^{p^k}.
\label{eq:d1p}\end{equation}

A higher derivation $(D_n)_{n\in\NN}$ is {\em iterative} \cite[p. 209]{MH} if
the following formula holds, for all integers $i,j\ge 0$:
\begin{equation}\label{eq:op1}
D_i\circ D_j=D_j
\circ D_i=\binomial{i+j}{i}D_{i+j}.
\end{equation}

It follows from (\ref{eq:op1}) and (\ref{eq:lucas}) with $r=1$ that
the operators $D_n$, $n\ge 0$, are completely determined by the
operators $D_{p^k}$ for $k\ge 0$. More precisely,
if $n=n_sp^s+\cdots+n_1p+n_0$ is the representation
of $n$ in base $p$ (so that $0\le n_i\le p-1$ for all $i$), then
the following formula holds:

\begin{equation}\label{eq:op3}
D_{n}=\frac{1}{n_s!\cdots n_0!}
D_{p^s}^{n_s}\circ\cdots\circ D_{p}^{n_1}\circ D_1^{n_0}.
\end{equation}

We now review some basic properties of hyperderivatives.
Let $f:\Omega\rightarrow C$ be a locally analytic function and
$z\in\Omega$ be a point. Recall that we have defined the hyperderivatives
$({\cal D}_nf)(z)$ of $f$ at $z$ by the formula
$$
f(z+\varepsilon)=\sum_{n\ge 0}({\cal D}_nf)(z)\varepsilon^n,\
\varepsilon\in C,\ \ |\varepsilon|\ \text{small}
$$
(this corresponds to Definition~2.3 of \cite{US}).

By Corollary 2.5
of \cite{US}, the family of operators $({\cal D}_n)_{n\in\NN}$ define a
higher derivation on the $C$-algebra of locally analytic functions, and it
is also well-known (and easy to check) that this higher derivation
is iterative.

Another important feature of the operators ${\cal D}_n$
is their good behaviour under holomorphy:
If $f$ is holomorphic on $\Omega$ and $n\in\NN$, then the function
${\cal D}_nf$ is also holomorphic on $\Omega$  by \cite[Theorem~3.1]{US},
and if moreover $f\in{\cal O}$, then it follows from \cite[Lemma~3.6]{US} that
${\cal  D}_nf\in {\cal O}$.

As explained in Section~\ref{intro}, it will be also convenient to
work with the following operators
\begin{equation}\label{operateur}
D_n:=\frac{1}{(-\overline{\pi})^n}{\cal D}_n\qquad (n\ge 0).
\end{equation}
The family $(D_n)_n$ is also an iterative higher derivation on ${\cal O}$, and
all the properties quoted above for the operators ${\cal D}_n$ remain
true for the operators $D_n$.

\subsection{Polynomial associated to a hyperderivative of a
quasi-modular form}\label{pol_associated}

We prove here Proposition~\ref{propo}.
This proposition will be a crucial tool in the subsequent
sections \ref{proof_th2}, \ref{h_hyperstable} and \ref{uptomodular}.

\begin{Proposition}\label{propo}
Let $f\in\qm_{w,m}^{\le\el}$ be a quasi-modular form of weight
$w$, type $m$ and depth $\le\el$, and let $P_f=\sum_{i=0}^{\el}f_iX^i$
be its associated polynomial. Then, for  
all $n\ge 0$,
${\cal D}_nf$ is an element of $\qm_{w+2n,m+n}^{\le\el+n}$,
with associated polynomial
$$
P_{{\cal D}_nf}=\sum_{j=0}^{n+\el}\Biggl[\sum_{r=0}^n
\binom{n+w+r-j-1}{r}{\cal D}_{n-r}f_{j-r}\Biggr]X^j,
$$
with the convention that $f_i=0$ if $i<0$ or $i>\el$.
\end{Proposition}

In the sequel, we will adopt the following notations. For
$\gamma=\begin{pmatrix}a & b\\ c & d\end{pmatrix}\in\Gamma_K$
and $z\in\Omega$, we write \[X_{\gamma}(z):=c/(cz+d),\quad
J_{\gamma}(z):=(cz+d)^2/(\det\gamma).\]
If $f:\Omega\rightarrow C$ is a holomorphic function, we will also
denote by $f\circ\gamma$ the function $z\mapsto f(\gamma(z))$.

The proof of Proposition~\ref{propo} requires three Lemmata.

\begin{Lemme}\label{binom}
Let $M,N,W$ be integers with $N\ge 0$. Then we have
\begin{eqnarray*}
 \sum_{i=0}^{N}(-1)^i\binom{M}{N-i}\binom{W+i-1}{i}&=&\binom{M-W}{N}\\
 & =& (-1)^N\binom{W+N-M-1}{N}.
\end{eqnarray*}
\end{Lemme}

\noindent\emph{Proof.} We have in $\FF_q[[X]]$:
$$ (1+X)^M=\sum_{i\ge 0}\binom{M}{i}X^i\quad\text{and}\quad
(1+X)^{-W}=\sum_{i\ge 0}(-1)^i\binom{W+i-1}{i}X^i,$$
so by multiplying we derive
$$ (1+X)^{M-W}=\sum_{N\ge 0}\left[\sum_{i=0}^N(-1)^i\binom{M}{N-i}
\binom{W+i-1}{i}\right] X^N.$$
Since on the other hand $(1+X)^{M-W}=\sum_{N\ge 0}\binom{M-W}{N}X^N$,
the result follows. \CVD

\begin{Lemme}\label{lemme5}
Let $f:\Omega\rightarrow C$ be a holomorphic function.
For all $n\ge 1$, $z\in\Omega$ and $\gamma\in\Gamma_K$,
we have
\begin{equation}\label{ccc9}
({\cal D}_n(f\circ\gamma))(z)=(-1)^nX_{\gamma}(z)^n\sum_{i=1}^n(-1)^i
\binom{n-1}{n-i}(X_{\gamma}(z)J_{\gamma}(z))^{-i}({\cal D}_if)(\gamma 
(z)).
\end{equation}
\end{Lemme}

\noindent\emph{Proof.}
Let $\varepsilon\in\Omega$ be sufficiently small.
We compute easily
$$ \gamma(z+\varepsilon)=\frac{\gamma(z)+\frac{a}{cz+d}\varepsilon}{1+
X_{\gamma}(z)\varepsilon}=\gamma(z)+H,$$
where
$$ H= J_{\gamma}(z)^{-1}\sum_{j\ge 1}(-1)^{j-1}X_{\gamma}(z)^{j-1}
\varepsilon^j.$$
Let us remark that $H$ tends to zero when $\varepsilon$ tends to zero.
We have for $i\ge 1$
$$
H^i=J_{\gamma}(z)^{-i}\sum_{n\ge i}(-1)^{n-i}X_{\gamma}(z)^{n-i}
\binom{n-1}{n-i}\epsilon^n,
$$
which gives
\begin{small}\begin{align*}
& f(\gamma(z+\varepsilon))=f(\gamma(z)+H)=f(\gamma(z))+\sum_{i\ge 1}
({\cal D}_if)(\gamma(z))H^i\\
& \phantom{(f\circ\gamma)(z+}
=f(\gamma(z))+\sum_{n\ge 1}(-1)^n\varepsilon^n\sum_{1\le i\le n}
(-1)^i\binom{n-1}{n-i}J_{\gamma}(z)^{-i}
X_{\gamma}(z)^{n-i}({\cal D}_if)(\gamma(z)).
\end{align*}\end{small}
Lemma~\ref{lemme5} follows at once.\CVD

\begin{Lemme}\label{lemme6}
Let $f$ be a non-zero element of $\qm_{w,m}^{\el}$,
and let $\gamma=\sqm{a}{b}{c}{d}\in\Gamma_K$.
Define $\varphi_{\gamma}:\Omega\rightarrow C$ by
$\varphi_{\gamma}(z):=(\det\gamma)^m(cz+d)^{-w}f(\gamma(z))$.
For $n\ge 0$ and $z\in\Omega$,
we have
\begin{equation}\label{ccc0}
({\cal D}_nf)(\gamma(z))=\frac{(cz+d)^{w+2n}}{(\det\gamma)^{m+n}}
\sum_{j=0}^{n}\binom{w+n-1}{n-j}X_{\gamma}(z)^{n-j}({\cal D}_j
\varphi_{\gamma})(z).
\end{equation}
\end{Lemme}

\noindent\emph{Proof.}
We argue by induction on $n$. For $n=0$ the Lemma is trivially true.
Let us now suppose that $n\ge 1$ and that the Lemma has been proved for
all integers $<n$. First of all, rewriting the formula (\ref{ccc9}) from
Lemma~\ref{lemme5}, we obtain
\begin{multline*}
({\cal D}_nf)(\gamma(z)) = J_{\gamma}(z)^n({\cal D}_n(f\circ\gamma)) 
(z)\\
{}+(-1)^{n-1}\sum_{i=1}^{n-1}(-1)^i
\binom{n-1}{n-i}X_{\gamma}(z)^{n-i} J_{\gamma}(z)^{n-i}({\cal D}_if) 
(\gamma(z)).
\end{multline*}
Next, using Lemma~\ref{lemme6} for $({\cal D}_if)(\gamma(z))$ (induction
hypothesis) and exchanging the priority of the summations thus  
obtained, we get
\begin{multline}
({\cal D}_nf)(\gamma(z))=J_{\gamma}(z)^n({\cal D}_n(f\circ\gamma))(z) 
\label{ccc10}\\
+(cz+d)^w(\det\gamma)^{-m}\sum_{j=0}^{n-1}\Biggl[
X_{\gamma}(z)^{n-j}J_{\gamma}(z)^n({\cal D}_j\varphi_{\gamma})(z)\\
{}\times
\sum_{i=j}^{n-1}(-1)^{n-1+i}\binom{n-1}{n-i}\binom{w+i-1}{i-j}\Biggr].
\end{multline}
In order to compute ${\cal D}_n(f\circ\gamma)$ in this formula,
we write $$ (f\circ\gamma)(z)=(\det\gamma)^{-m}(cz+d)^w\varphi_ 
{\gamma}(z)$$
and apply Leibniz's rule (\ref{eq:op2}). We find
\begin{align*}
({\cal D}_n(f\circ\gamma))(z)&={\cal D}_n\bigl((\det\gamma)^{-m}
(cz+d)^w\varphi_{\gamma}(z)\bigr)\\
& = (\det\gamma)^{-m}\sum_{j=0}^n{\cal D}_{n-j}\bigl((cz+d)^w\bigr)
({\cal D}_j\varphi_{\gamma})(z)\\
&=(\det\gamma)^{-m}(cz+d)^w\sum_{j=0}^n\binom{w}{n-j}X_{\gamma}(z)^{n-j}
({\cal D}_j\varphi_{\gamma})(z),
\end{align*}
the last line coming from the following formula, which is easy to check:
$$
{\cal D}_i\bigl((cz+d)^w\bigr)=\binom{w}{i}X_{\gamma}(z)^i(cz+d)^w\quad
(i\ge 0).
$$
Substituting in (\ref{ccc10}) the expression obtained for
${\cal D}_n(f\circ\gamma)$, we derive
\begin{align*}
({\cal D}_nf)(\gamma(z)) & = \frac{(cz+d)^w}{(\det\gamma)^m}J_{\gamma} 
(z)^n
\sum_{j=0}^nX_{\gamma}(z)^{n-j}({\cal D}_j\varphi_{\gamma})(z)
\Biggr[\binom{w}{n-j}\\
&\phantom{= \frac{(cz+d)^w}{(\det\gamma)^m}}
{}+(-1)^{n-1}\sum_{i=0}^{n-1-j}(-1)^{i+j}
\binom{n-1}{n-i-j}\binom{w+i+j-1}{i}\Biggl],
\end{align*}
which gives, using Lemma~\ref{binom} with $N=n-j$, $M=n-1$ and $W=w+j$,
$$
({\cal D}_nf)(\gamma(z)) = \frac{(cz+d)^{w+2n}}{(\det\gamma)^{m+n}}
\sum_{j=0}^n\binom{w+n-1}{n-j}X_{\gamma}(z)^{n-j}({\cal D}_j
\varphi_{\gamma})(z).
$$
This is the desired formula for $n$.\CVD

\noindent\emph{Proof of Proposition~\ref{propo}.} Let $\gamma$ be an element
of $\Gamma_K$ and let $\varphi_{\gamma}$ be defined as in Lemma~\ref{lemme6}.
Then we have, for all $z\in\Omega$:
\begin{equation}\label{ccc1}
\varphi_{\gamma}(z)=\sum_{i=0}^{\el}f_i(z)X_{\gamma}(z)^i.
\end{equation}
The following formula holds for all $i,r\ge 0$, and is easy to check:
\begin{equation}
{\cal D}_r(X_{\gamma}^i)=(-1)^r\binom{i+r-1}{r}X_{\gamma}^{i+r}.
\label{eq:der_X_gamma}\end{equation}
Thus, Leibniz's formula (\ref{eq:op2}) gives, for all $j$ such that
$0\le j\le n$,
\begin{equation}\label{ccc2}
{\cal D}_j(f_iX_{\gamma}^i)=\sum_{r=0}^{j}(-1)^{j-r}
\binom{i+j-r-1}{j-r}X_{\gamma}^{i+j-r}{\cal D}_rf_i.
\end{equation}
Using now successively (\ref{ccc0}), (\ref{ccc1}) and (\ref{ccc2}), we obtain
\begin{align*}
\frac{(\det\gamma)^{m+n}}{(cz+d)^{w+2n}} ({\cal D}_n&f)(\gamma(z)) =
\sum_{j=0}^{n}\binom{k+n-1}{n-j}X_{\gamma}(z)^{n-j}
\sum_{i=0}^{\el}({\cal D}_j(f_iX_{\gamma}^i))(z)\\
= \sum_{\genfrac{}{}{0pt}{1}{0\le j\le n}{
\genfrac{}{}{0pt}{1}{0\le r\le j}{0\le i\le\el}}} &
(-1)^{j-r}\binom{k+n-1}{n-j}
\binom{i+j-r-1}{j-r}X_{\gamma}(z)^{n+i-r}({\cal D}_rf_i)(z)\\
= \sum_{\genfrac{}{}{0pt}{1}{0\le u\le n+\el}{
\genfrac{}{}{0pt}{1}{0\le r\le n}{r\le j\le n}}} &
(-1)^{j-r}\binom{w+n-1}{n-j}
\binom{u-n+j-1}{j-r}X_{\gamma}(z)^{u}({\cal D}_rf_{u-n+r})(z).
\end{align*}
This formula shows already that
$(\det\gamma)^{m+n}(cz+d)^{-(k+2n)}({\cal D}_nf)(\gamma(z))$ is
a polynomial in $X_{\gamma}(z)$ of degree $\le n+\el$ with
coefficients in $\cal O$, whence
${\cal D}_nf\in\qm_{w+2n,m+n}^{\le \el+n}$. Finally, we get the desired
formula for $P_{{\cal D}_nf}$ by making the change of summation index
$r'=n-r$ in the last formula, and then using Lemma~\ref{binom} with
$N=r'$, $M=w+n-1$, $W=u-r'$.\CVD

\subsection{Proof of Theorem~\ref{theo2}}\label{proof_th2}

As a first Corollary of Proposition~\ref{propo}, we will now derive
Theorem~\ref{theo2}. We will need the following Lemma, which will also be
useful later.

\begin{Lemme}\label{lemme7}
Let $f:\Omega\rightarrow C$ be an element of the ring $\cal O$
having the $t$-expansion
\begin{equation}
f(z) = \sum_{n\ge 0}a_n t^n, \label{texpansion}
\end{equation}
where $t=t(z)$.
Then, for all integers $i\ge 1$, one has
$$
(D_if)(z)=\sum_{n\ge 2}b_{i,n}t^n,
$$
where
$$
b_{i,n}=\sum_{r=1}^{n-1}(-1)^{i+r}\binom{n-1}{r}
\Bigl(
\sum_{\genfrac{}{}{0pt}{1}{i_1,\ldots,i_r\ge 0}{q^{i_1}+\cdots+q^{i_r} 
=i}}
\frac{1}{d_{i_1}\cdots d_{i_r}}\Bigr)a_{n-r}.
$$
\end{Lemme}

\noindent{\it Proof.}
Fix $z\in\Omega$ with $|z|_i$ sufficiently large, so that
$f$ has a $t$-expansion (\ref{texpansion}), and let $\varepsilon$
be any element of the field $K_{\infty}$. Then
$|z+\varepsilon|_i=|z|_i$, and thus
\begin{equation}\label{ccc11}
   f(z+\varepsilon)=\sum_{n\ge 0}a_nt(z+\varepsilon)^n.
\end{equation}
Now, putting as usual $t:=t(z)$ and using the formula (\ref{exponential}),
we have, for $|\varepsilon|$ sufficiently small:
\begin{align}
t(z+\varepsilon)^n=\frac{1}{\bigl(e_{{\mathcal C}}(\overline{\pi} z)+
e_{{\mathcal C}}(\overline{\pi}\varepsilon)\bigr)^n} & = \frac{t^n}{\bigl(1
+e_{{\mathcal C}}(\overline{\pi}\varepsilon)t\bigr)^n}\nonumber\\
& = \sum_{r\ge 0}(-1)^r\binom{r+n-1}{r}(e_{{\mathcal C}}(\overline{\pi}
\varepsilon))^rt^{r+n}
\nonumber\\
& = \sum_{r\ge 0}(-1)^r\binom{r+n-1}{r}t^{r+n}\Bigl(\sum_{i\ge 0}
\frac{\overline{\pi}^{q^i}\varepsilon^{q^i}}{d_i}\Bigr)^r\nonumber\\
& = \sum_{r\ge 0}(-1)^r\binom{r+n-1}{r}t^{r+n}\sum_{i\ge 0}
\alpha_{r,i}\overline{\pi}^i\varepsilon^i\label{ccc12},
\end{align}
where
$$
\alpha_{r,i} =
\sum_{\genfrac{}{}{0pt}{1}{i_1,\ldots,i_r\ge 0}{q^{i_1}+\cdots+q^{i_r}=i}}
\frac{1}{d_{i_1}\cdots d_{i_r}}\qquad\text{if\ } r\ge 1\ \text{and\ }  
i\ge 0,
$$
and where
$$
\alpha_{0,i} =\left\{
\begin{array}{cl}
1 & \text{for\ } i=0\\
0 & \text{for\ } i\not=0.
\end{array}\right.
$$
Substituting (\ref{ccc12}) in (\ref{ccc11}) and exchanging the order of
the summations, we get
\begin{equation*}
f(z+\varepsilon)=\sum_{i\ge 0}\Biggl[\sum_{r\ge 0}\sum_{n\ge 0}(-1)^ra_n
\alpha_{r,i}\binom{r+n-1}{r}t^{r+n}\Biggr]\overline{\pi}^i\varepsilon^i.
\end{equation*}
Since on the other hand we have (by virtue of (\ref{operateur}))
\begin{equation*}
f(z+\varepsilon)=\sum_{i\ge 0}(D_if)(z)(-1)^i\overline{\pi}^i 
\varepsilon^i,
\end{equation*}
we obtain
\begin{equation}\label{ccc13}
(D_if)(z)=\sum_{r\ge 0}\sum_{n\ge 0}(-1)^{r+i}\binom{r+n-1}{r}\alpha_ 
{r,i}a_n
t^{r+n}\quad (i\ge 0).
\end{equation}
Let now $i\ge 1$. For $r=0$ we have $\alpha_{r,i}=0$, and for
$r\ge 1$ and $n=0$ we have $\binom{r+n-1}{r}=0$. Thus, when $i\ge 1$,
we can take $r\ge 1$ and $n\ge 1$ in the formula
(\ref{ccc13}). This gives the Lemma. \CVD

We are now in position to prove Theorem~\ref{theo2}.

\medskip

\noindent\emph{Proof of Theorem~\ref{theo2}.} By Corollary~\ref{coroll},
it suffices to prove that $D_n(\qm_{w,m}^{\le\el}(K))\subset
\qm_{w+2n,m+n}^{\le\el+n}(K)$.
But if $f\in\qm_{w,m}^{\le\el}(K)$, then the function $D_nf$ is an  
element
of $\qm_{w+2n,m+n}^{\le\el+n}$ by Proposition~\ref{propo}, and the
coefficients of its $t$-expansion all lie in $K$ after Lemma~\ref 
{lemme7}.
Hence $f\in\qm_{w+2n,m+n}^{\le\el+n}(K)$ as claimed. \CVD

\subsection{The principal ideal $(h)$ is hyperdifferential.}\label{h_hyperstable}

As a second application of Proposition \ref{propo}, we prove here
the following Proposition, which will be very useful in this text.

\begin{Proposition}\label{hstable}
The ideal $(h)$ is hyperdifferential.
\end{Proposition}

\noindent\emph{Proof.} Applying Proposition~\ref{propo} with
$f=h\in M_{q+1,1}$ we see
that, for all $n\ge 0$, the polynomial $P_{{\cal D}_nh}$ associated to
${\cal D}_nh$ is:
\begin{equation}\label{ccc14}
P_{{\cal D}_nh}=\sum_{j=0}^n\binom{n+q}{j}({\cal D}_{n-j}h)X^j.
\end{equation}
Now, the function ${\cal D}_{n-j}h$ vanishes at infinity if $0\le j<n$
after Lemma~\ref{lemme7}, and of course the same holds for $j=n$  
because $\nu_\infty(h)=1$. Since 
$h$ does
not vanish on $\Omega$ (by formula (5.14) of \cite{Ge}),
the  
functions $f_j:={\cal D}_{n-j}h/h$
($0\le j\le n$) are all holomorphic on $\Omega$ and at infinity;  
consequently, they belong to $\cal O$. We deduce from
(\ref{ccc14}) that for all $n\ge 0$, $z\in\Omega$ and
$\gamma\in \Gamma_K$:
$$ \frac{{\cal D}_nh}{h}(\gamma(z))=(cz+d)^{2n}(\det\gamma)^{-n}
\sum_{j=0}^n\binom{n+q}{j}f_j(z)X_{\gamma}(z)^j, $$
with $f_j\in\cal O$. In other words, the function ${\cal D}_nh/h$
is a quasi-modular form (of weight $2n$, type $n$ and depth $\le n$),
and therefore belongs to the polynomial ring $C[E,g,h]$ by virtue of  
Theorem~\ref{theo1}.\CVD

\subsection{Hyperderivatives of $E,g,h$.}\label{uptomodular}

We end this Section~\ref{section_hyp} by
establishing Proposition~\ref{derivees_mod_modulaires} below, which is another
consequence of Proposition~\ref{propo}. Roughly speaking,
this Proposition determines, up to a modular form, the functions
$D_{p^i}E$, $D_{p^i}g$ and $D_{p^i}h$ for all integers $i\ge 0$.
As we will see in Section~\ref{calculs}, this result, together with some
technical computations, yields a practical way to compute explicitly
$D_{p^i}E$, $D_{p^i}g$ and $D_{p^i}h$ when $p^i\le q^2$. Proposition
\ref{derivees_mod_modulaires} will
also be needed in the proof of Lemma~\ref{lemma:munu}.
\begin{Proposition}\label{derivees_mod_modulaires}
\begin{enumerate}
\item[\emph{(i)}] For every integer $i\ge 0$, we have
$$D_{p^i}E-E^{p^i+1}\in M_{2p^i+2,p^i+1}.$$
\item[\emph{(ii)}] For every integer $i\ge 1$, we have
$$\left\{\begin{array}{l}
D_{p^i}g=0\ \ \mbox{\text if}\ \  p^i<q\\
D_{p^i}g-E^{p^i}g\in M_{2p^i+q-1,p^i}\ \  \mbox{\text if}\ \ p^i\ge q.
\end{array}\right.$$
\item[\emph{(iii)}] For every integer $i\ge 0$, we have
$$ D_{p^i}h-E^{p^i}h-E^q(D_{p^i-q}h)\in M_{2p^i+q+1,p^i+1},$$
where we have used the convention $D_nh=0$ when $n<0$.
\end{enumerate}
\end{Proposition}

To prove Proposition~\ref{derivees_mod_modulaires}, we will need
three Lemmata.

\begin{Lemme}\label{fonction_E}
For every integer $i\ge 0$, we have $D_{p^i-1}E=E^{p^i}$.
\end{Lemme}

\noindent\emph{Proof.} To prove this Lemma, we use the
following formula
(see \cite[page 686]{Ge}), valid for every $z\in\Omega$:
$$ E(z)=\frac{1}{\overline{\pi}}
\sum_{\genfrac{}{}{0pt}{1}{a\in A}{a\,\text{monic}}}\sum_{b\in A}
\frac{a}{az+b}.$$
Since the series are here uniformly convergent (in the
above order) on every open ball contained in $\Omega$, we find,
by \cite[Proposition 2.6]{US} and by using the formula
$$ {\cal D}_n\bigl(\frac{a}{az+b}\bigr)=(-1)^n\bigl(
\frac{a}{az+b}\bigr)^{n+1}\qquad (n\ge 0)$$ which follows from
(\ref{eq:der_X_gamma}),
\begin{equation*}
({\cal D}_{p^i-1}E)(z) = \frac{1}{\overline{\pi}}
\sum_{a\ \text{monic}}\sum_{b\in A}(-1)^{p^i-1}\bigl(
\frac{a}{az+b}\bigr)^{p^i}=(-\overline{\pi})^{p^i-1} E(z)^{p^i}.
\end{equation*}
Thus, we get $D_{p^i-1}E=E^{p^i}$ as annouced.\CVD

\begin{Lemme}\label{fonction_g}
We have $D_ng=0$ for every integer $n$ such that $2\le n<q$.
\end{Lemme}

\noindent\emph{Proof.} By \cite[formula (6.8)]{Ge} we have, for all
$z\in\Omega$:
$$
g(z)=\overline{\pi}^{1-q}(T^q-T)
\sum_{\genfrac{}{}{0pt}{1}{(a,b)\in A^2}{(a,b)\not=0}}\frac{1}{(az+b)^{q-1}}.
$$
Using the formula (\ref{eq:der_X_gamma}), we deduce that for all
$n\ge 0$ the following holds:
$$
({\cal D}_ng)(z)=\overline{\pi}^{1-q}(T^q-T)(-1)^n\binom{n+q-2}{n}
\sum_{\genfrac{}{}{0pt}{1}{(a,b)\in A^2}{(a,b)\not=0}}
\frac{a^n}{(az+b)^{n+q-1}}.
$$
Since now for $2\le n\le q-1$ we have in $\FF_p$ (for instance by
(\ref{eq:lucas}) with $p^r=q$)
$\binom{n+q-2}{n}=0$, we get the Lemma. \CVD

\begin{Lemme}\label{fonction_h}
For every integer $i\ge 0$ such that $p^i\ge q$, we have
$D_{p^i-q}h\in M_{2p^i-q+1,p^i}$.
\end{Lemme}

\noindent\emph{Proof.} Applying Proposition~\ref{propo} we get
$$
P_{{\cal D}_{p^i-q}h}=\sum_{j=0}^{p^i-q}\binom{p^i}{j}({\cal D}_{p^i- 
q-j}h)X^j.
$$
But in $\FF_p$ we have (see (\ref{eq:lucas}))
\begin{equation*}
\binom{p^i}{j}=\left\{\begin{array}{ccl}
1 & \text{if} & j=0\\
0 & \text{if} & 0<j<p^i.
\end{array}\right.
\end{equation*}
Therefore, we obtain $P_{{\cal D}_{p^i-q}h}={\cal D}_{p^i-q}h$.
Thus, $\deg (P_{{\cal D}_{p^i-q}h})=0$ and it follows that
${\cal D}_{p^i-q}h$ is a modular form, which proves the Lemma. \CVD

\noindent\emph{Proof of Proposition~\ref{derivees_mod_modulaires}.}
(i) For $i=0$ the assertion (i) is clearly true by
the relations (\ref{system}), so we may assume $i\ge 1$.
Since $P_E=-\overline{\pi}^{-1}X+E$, we first note from
Proposition~\ref{propo} that the polynomial $P_{{\cal D}_{p^i}E}$
associated to ${\cal D}_{p^i}E$ equals
$$ P_{{\cal D}_{p^i}E} = -\frac{1}{\overline{\pi}}X^{p^i+1}+
\sum_{j=0}^{p^i}\binom{p^i+1}{j}({\cal D}_{n-j}E)X^j.$$
But in $\FF_p$ we have
\begin{equation*}
\binom{p^i+1}{j}=\left\{\begin{array}{ccl}
1 & \text{if} & j=0, j=1\ \text{or}\ j=p^i\\
0 & \text{if} & 1<j<p^i.
\end{array}\right.
\end{equation*}
Therefore, by Lemma~\ref{fonction_E},
\begin{align*}
P_{{\cal D}_{p^i}E} &= -\frac{1}{\overline{\pi}}X^{p^i+1}+
EX^{p^i}+({\cal D}_{p^i-1}E)X+{\cal D}_{p^i}E\\
& = -\frac{1}{\overline{\pi}}X^{p^i+1}+EX^{p^i}+
(-\overline{\pi})^{p^i-1}E^{p^i}X+{\cal D}_{p^i}E.
\end{align*}
Arguing now as in the proof of Lemma~\ref{lemme4}, if we define
$f=(-\overline{\pi})^{-p^i}{\cal D}_{p^i}E-E^{p^i+1}$, we get:
\begin{align*}
P_f & = (-\overline{\pi})^{-p^i}P_{{\cal D}_{p^i}E}-
(-\frac{1}{\overline{\pi}}X+E)^{p^i+1}\\
& = (-\overline{\pi})^{-p^i}P_{{\cal D}_{p^i}E}-
((-\overline{\pi})^{-p^i}X^{p^i}+E^{p^i})(-\frac{1}{\overline{\pi}}X 
+E)=f.
\end{align*}
Since $\deg P_f=0$, it follows that $f$ is a modular form, which
proves (i).

(ii) Let $i\ge 1$ be an integer. If $p^i<q$, then the assertion (ii) is
true by Lemma~\ref{fonction_g}, so we may now assume $p^i\ge q$.
Applying proposition~\ref{propo}, we have
\begin{equation}
P_{{\cal D}_{p^i}g}=\sum_{j=0}^{p^i}\binom{p^i+q-2}{j}({\cal D}_{p^i-j}g)X^j.
\label{associated_gg}
\end{equation}
For every integer $j$ such that $0\le j\le p^i$, write $p^i-j=\alpha q+
\beta$ with $\alpha\ge 0$ and $0\le\beta\le q-1$.
If $\beta\ge 2$, then ${\cal D}_{\beta}g=0$ by Lemma~\ref{fonction_g},
and thus
$$ {\cal D}_{p^i-j}g = {\cal D}_{\alpha q+\beta}g =
\binom{\alpha q+\beta}{\beta}{\cal D}_{\alpha q+\beta}g=
{\cal D}_{\alpha q}({\cal D}_{\beta}g)=0.$$
If $\beta=1$, one has
$$ \binom{p^i+q-2}{j}=\binom{p^i+q-2}{\alpha q+q-1}=0, $$
and if $\beta=0$, one finds easily (recall that $q$ is a power of $p$)
$$ \binom{p^i+q-2}{j}=\binom{p^i+q-2}{\alpha q+q-2}=\left\{
\begin{array}{ccl}
1 & \text{if} & \alpha=0\text{\ or\ }\alpha=p^i/q\\
0 & \text{if} & 0<\alpha <p^i/q.
\end{array}\right.$$
Putting these remarks together, (\ref{associated_gg}) leads to
$$ P_{{\cal D}_{p^i}g}=gX^{p^i}+{\cal D}_{p^i}g.$$
Thus, by the same argument as before,
$D_{p^i}g-E^{p^i}g$ is a modular form, and so is an element
of $M_{2p^i+q-1,p^i}$.

(iii) For every integer $i\ge 0$ we have, by Proposition~\ref{propo},
\begin{equation}\label{associated_h}
P_{{\cal D}_{p^i}h}=\sum_{j=0}^{p^i}\binom{p^i+q}{j}({\cal D}_{p^i-j}h)X^j.
\end{equation}
In the case $p^i<q$ one easily computes
$$ \binom{p^i+q}{j}=\left\{\begin{array}{ccl}
0 & \text{if} & 0<j<p^i\\
1 & \text{if} & j=0 \text{\ or\ }j=p^i.
\end{array}\right.$$
If $p^i=q$, then
\begin{equation*}
\binom{p^i+q}{j}=\left\{\begin{array}{ccl}
1 & \text{if} & j=0\\
2 & \text{if} & j=p^i\\
0 & \text{if} & 0<j<p^i,
\end{array}\right.
\end{equation*}
and if $p^i>q$, we find
\begin{equation*}
\binom{p^i+q}{j}=\left\{\begin{array}{ccl}
1 & \text{if} & j\in\{0,q,p^i\}\\
0 & \text{if} & 0<j<p^i, j\not=q.
\end{array}\right.
\end{equation*}
Hence, we deduce in the three cases from (\ref{associated_h}):
$$ P_{{\cal D}_{p^i}h}={\cal D}_{p^i}h+({\cal D}_{p^i-q}h)X^q+hX^ 
{p^i}.$$
If we now set
$f=(-\overline{\pi})^{-p^i}{\cal D}_{p^i}h-E^{p^i}h$, then
we have
\begin{align*}
P_f & = (-\overline{\pi})^{-p^i}P_{{\cal D}_{p^i}h}-
(-\frac{1}{\overline{\pi}}X+E)^{p^i}h\\
& = (-\overline{\pi})^{-p^i}
\bigl({\cal D}_{p^i}h+({\cal D}_{p^i-q}h)X^q+hX^{p^i}\bigr)
-(-\overline{\pi})^{-p^i}hX^{p^i}-E^{p^i}h\\
& = f + (-\overline{\pi})^{-p^i}({\cal D}_{p^i-q}h)X^q.
\end{align*}
Since the function $D_{p^i-q}h$ is a modular form by
Lemma~\ref{fonction_h}, it follows that $f-E^q(D_{p^i-q}h)$ is also
modular.
This means that $$D_{p^i}h-E^{p^i}h-E^q(D_{p^i-q}h)\in M_{2p^i+q+1,p^i 
+1}$$
and part (iii) of Proposition~\ref{derivees_mod_modulaires} is proved.
\CVD

\noindent\textbf{Remark.} Proposition \ref{derivees_mod_modulaires} may be
applied to give another proof of Proposition \ref{hstable}.
Due to the importance of the latter result in this text,
we give here all the details.

\medskip

\noindent\emph{Alternative proof of Proposition \ref{hstable}.}
By Proposition \ref{derivees_mod_modulaires} Part (iii), the function
$D_{p^i}h-E^{p^i}h-E^q(D_{p^i-q}h)$ is a modular form for every $i\geq 0$.
By Lemma \ref{lemme7}, this modular form vanishes at infinity, hence belongs
to the ideal $(h)$. From this we deduce:
\[\begin{array}{rclr}D_{p^i}h &\in & (h)& \mbox{ for } p^i<q\\ 
D_{p^i}h-E^qD_{p^i-q}h &\in & (h)& \mbox{ for } p^i\geq q. \end{array}\]

By induction on $i\geq 0$, we now show that $D_{n}h\in(h)$ for
all $0\leq n<p^i$. The first property above, together with the formulas
(\ref{eq:op2}) and (\ref{eq:op3}), easily imply that $D_nh\in (h)$
for all $n<q$, so the assertion is true for $i$
such that $1\leq p^i\le q$. Supposing now that the property holds for an
integer $i$ such that $p^i\geq q$, we have in particular
that $D_{p^i-q}h\in(h)$, which implies $D_{p^i}h\in(h)$.
By (\ref{eq:op2}), (\ref{eq:op3}) and the induction hypothesis, we then get
$D_nh\in(h)$ for all $0\leq n<p^{i+1}$, {\em i.e.} the assertion for $i+1$.\CVD

\section{Explicit formulas for $D_i$ when $1\le i\le q^2$\label{calculs}}

In this (more technical) Section, we compute explicitly $D_nE$, $D_ng$ and
$D_nh$ as polynomials in $E$, $g$, $h$, when $n\le q^2$ is a power of $p$.
The formulas obtained allow to compute (at least in principle) $D_nf$ for
any $f\in C[E,g,h]$ and any $n$ such that $0\le n\le q^2$, thanks to
Leibniz's rule (\ref{eq:op2}) and the formula (\ref{eq:op3}).

We summarise the main results in the following Theorem. Here, $e$ is
the integer such that $q=p^e$.

\begin{Theorem}\label{theo_calculs}
Let $i,n$ be integers. The following formulas hold.
\begin{itemize}
         \item[\emph{(i)}] $D_nE=E^{n+1}$ if $0\le n<q$.
         \item[\emph{(ii)}] $D_1g=-(Eg+h)$ and $D_ng=0$ if $2\le n<q$.
         \item[\emph{(iii)}] $D_nh=E^nh$ if $0\le n<q$.
         \item[\emph{(iv)}] $D_{p^i}E=E^{p^i+1}+\frac{1}{d_1^{p^{i-e}}}
          g^{p^{i-e}-1}h^{p^{i-e}+1}$ if $q\le p^i<q^2$.
         \item[\emph{(v)}] $D_{p^i}g=E^{p^i}g$ if $q\le p^i<q^2$.
         \item[\emph{(vi)}] $D_{p^i}h=E^{p^i}h+\frac{1}{d_1^{p^{i-e}-1}}E^q
         g^{p^{i-e}-1}h^{p^{i-e}}-\frac{1}{d_1^{p^{i-e}}}g^{p^{i-e}}
         h^{p^{i-e}+1}$ if $q\le p^i<q^2$.
         \item[\emph{(vii)}] $D_{q^2}E=E^{q^2+1}+\frac{1}{d_1^q}g^{q-1}h^{q+1}
         +\frac{1}{d_2}g^{2q}h^{2}$.
         \item[\emph{(viii)}] $D_{q^2}g=E^{q^2}g-\frac{d_1}{d_2}g^{q+1}h^q
         +(\frac{1}{d_1^{q-1}}-\frac{d_1^2}{d_2})h^{2q-1}$.
         \item[\emph{(ix)}] $D_{q^2}h=E^{q^2}h+\frac{1}{d_1^{q-1}}
         E^qg^{q-1}h^q-\frac{1}{d_2}g^{2q+1}h^2-(\frac{d_1}{d_2}+
         \frac{1}{d_1^q})g^qh^{q+1}$.
\end{itemize}
\end{Theorem}

The proof of Theorem~\ref{theo_calculs} requires several technical  
Lemmata.

\begin{Lemme}\label{developp}
Let $i,j$ be integers such that $q< p^i<q^2$ and $q<p^j\le q^2$.
We have the following $t$-expansions.
\begin{enumerate}
\item[\emph{(i)}] $E=t+t^{q^2-2q+2}+\cdots$
\item[\emph{(ii)}] $g=1-d_1t^{q-1}-d_1t^{q^3-2q^2+2q-1}+\cdots$
\item[\emph{(iii)}] $h=-t-t^{q^2-2q+2}+\cdots$
\item[\emph{(iv)}] $D_qE=\frac{1}{d_1}t^2+t^{q+1}+\cdots$
\item[\emph{(v)}] $D_qg=t^q+\cdots$
\item[\emph{(vi)}] $D_qh=-\frac{1}{d_1}t^2-t^{q+1}+\cdots$
\item[\emph{(vii)}] $D_{p^i}E=\frac{1}{d_1^{p^{i-e}+1}}t^{p^{i-e}+1}+ 
\cdots$
\item[\emph{(viii)}] $D_{p^i}g=O(t^{p^{i-e}+1})$
\item[\emph{(ix)}] $D_{p^i}h=-\frac{1}{d_1^{p^{i-e}+1}}t^{p^{i-e}+1}+ 
\cdots$
\item[\emph{(x)}] $D_{p^j-q}h=-\frac{1}{d_1^{p^{j-e}-1}}t^{p^{j-e}}+ 
\cdots$
\item[\emph{(xi)}] $D_{q^2}E=\frac{1}{d_2}t^2+\frac{1}{d_1^q}t^{q+1}+ 
\cdots$
\item[\emph{(xii)}] $D_{q^2}g=\frac{d_1}{d_2}t^q-\frac{1}{d_1^{q-1}}t^ 
{2q-1}+\cdots$
\item[\emph{(xiii)}] $D_{q^2}h=-\frac{1}{d_2}t^2-\frac{1}{d_1^q}t^{q 
+1}+\cdots$
\end{enumerate}
\end{Lemme}

\noindent \emph{Proof.} The formulas (i), (ii) and (iii) follow
from \cite{Ge}, Corollaries (10.5), (10.11) and (10.4).
To determine the other $t$-expansions, we apply Lemma~\ref{lemme7}.
Define, for $r\ge 1$ and $i\ge 1$,
$$
\alpha_{r,i} =
\sum_{\genfrac{}{}{0pt}{1}{i_1,\ldots,i_r\ge 0}{q^{i_1}+\cdots+q^{i_r} 
=i}}
\frac{1}{d_{i_1}\cdots d_{i_r}}.
$$
One computes painlessly, for all integers $i,j$ such that $q<p^i<q^2$
and $q<p^j\le q^2 $,
\begin{equation*}
\alpha_{r,q}=\left\{
\begin{array}{ccl}
    1/d_1 & \text{if} & r=1\\
    0     & \text{if} & 2\le r\le q-1 \\
    1/d_1^q & \text{if} & r=q,
\end{array}\right.
\end{equation*}
\begin{equation*}
\alpha_{r,p^i}=\left\{
\begin{array}{ccl}
    0 & \text{if} & 1\le r < p^{i-e}\\
    1/d_1^{p^{i-e}} & \text{if} & r=p^{i-e},
\end{array}\right.
\end{equation*}
\begin{equation*}
\alpha_{r,p^j-q}=\left\{
\begin{array}{ccl}
    0 & \text{if} & 1\le r < p^{j-e}-1\\
    1/d_1^{p^{j-e}-1} & \text{if} & r=p^{j-e}-1.
\end{array}\right.
\end{equation*}
and
\begin{equation*}
\alpha_{r,q^2}=\left\{
\begin{array}{ccl}
    1/d_2 & \text{if} & r=1\\
    0     & \text{if} & 2\le r\le q-1 \\
    1/d_1^q & \text{if} & r=q\\
    0     & \text{if} & q < r < 2q-1\\
    1/d_1^{q-1} & \text{if} & r=2q-1.\\
\end{array}\right.
\end{equation*}
Thanks to these expressions and the expansions (i), (ii), (iii)
already obtained, Lemma~\ref{lemme7} now yields Lemma~\ref{developp}
without difficulty.\CVD

In the next Lemma, if $v_1,\ldots,v_n$ are elements of a vector space,
we denote by $\langle v_1,\ldots,v_n\rangle $ the subspace they  
generate.

\begin{Lemme}\label{espacesM}
Let $i,n$ be integers. The following holds:
\begin{enumerate}
\item[\emph{(i)}] If $p\le p^i < q$, then
$M_{2p^i+2,p^i+1}=\left\{\begin{array}{lcl}
                   \{0\} & \text{if} & (p,q,i)\not=(2,4,1) \\
                   \langle g^2\rangle  & \text{if} & (p,q,i)=(2,4,1).
\end{array}\right.$
\item[\emph{(ii)}] If $q\le p^i < q^2$, then
$$M_{2p^i+2,p^i+1}=\left\{\begin{array}{ll}           
           \langle g^6,gh^3\rangle  & \text{if} \ (p,q,i)=(2,4,3)\\
           \langle g^4,h^2\rangle  & \text{if} \  (p,q,i)=(3,3,1)\\
           \langle g^6,g^3h,h^2\rangle  & \text{if} \ (p,q,i)=(2,2,1)\\
\langle g^{p^{i-e}-1}h^{p^{i-e}+1}\rangle  & \text{otherwise.}
\end{array}\right.$$
\item[\emph{(iii)}]
$M_{2q^2+2,q^2+1}=\left\{\begin{array}{lcl}
                   \langle g^{q-1}h^{q+1},g^{2q}h^2\rangle  & \text{if}
& q\ge 4 \\
                   \langle g^2h^4,g^6h^2,g^{10}\rangle  & \text{if} & q=3 \\
                   \langle gh^3,g^4h^2,g^7h,g^{10}\rangle  & \text{if} & q=2.
\end{array}\right.$
\item[\emph{(iv)}] If $q\le p^i < q^2$, then
$$M_{2p^i+q-1,p^i}=\left\{\begin{array}{lcl}
                   \langle g^{p^{i-e}+1}h^{p^{i-e}}\rangle & \text{if}
& (p,q,i)\not=(2,2,1) \\
                   \langle g^5,g^2h\rangle  & \text{if} & (p,q,i)=(2,2,1).
\end{array}\right.$$
\item[\emph{(v)}]
$M_{2q^2+q-1,q^2}=\langle g^{2q+2}h,g^{q+1}h^q,h^{2q-1}\rangle $.
\item[\emph{(vi)}] If $p\le p^i < q$, then
$M_{2p^i+q+1,p^i+1}=\left\{\begin{array}{lcl}
                   \{0\} & \text{if} & (p,q,i)\not=(2,4,1) \\
                   \langle g^3\rangle  & \text{if} & (p,q,i)=(2,4,1) \\
\end{array}\right.$
\item[\emph{(vii)}] If $q\le p^i < q^2$, then
$$M_{2p^i+q+1,p^i+1}=\left\{\begin{array}{ll}
           \langle g^7,g^2h^3\rangle  & \text{if}\ (p,q,i)=(2,4,3)\\
           \langle g^5,gh^2\rangle  & \text{if}\ (p,q,i)=(3,3,1)\\
           \langle g^7,g^4h,gh^2\rangle  & \text{if}\ (p,q,i)=(2,2,1)\\
\langle g^{p^{i-e}}h^{p^{i-e}+1}\rangle  & \text{otherwise.}
\end{array}\right.$$
\item[\emph{(viii)}]
$M_{2q^2+q+1,q^2+1}=\left\{\begin{array}{lcl}
                   \langle g^{q}h^{q+1},g^{2q+1}h^2\rangle  & \text 
{if} & q\ge 4 \\
                   \langle g^3h^4,g^7h^2,g^{11}\rangle  & \text{if} &  
q=3 \\
                   \langle g^2h^3,g^5h^2,g^8h,g^{11}\rangle  & \text 
{if} & q=2.
\end{array}\right.$
\item[\emph{(ix)}]
$M_{2p^i-q+1,p^i}=\left\{\begin{array}{lcl}
                   \langle g^{p^{i-e}-1}h^{p^{i-e}}\rangle  & \text 
{if} & q< p^i < q^2\\
                   \langle g^{q-1}h^q,g^{2q}h\rangle  & \text{if} &  
p^i=q^2.
\end{array}\right.$
\end{enumerate}
\end{Lemme}

\noindent \emph{Proof.} We know that a basis of $M_{w,m}$
consists of monomials $g^{\alpha}h^{\beta}$, where $\alpha$ and $\beta$
are non-negative integers satisfying
\begin{equation*}
\left\{\begin{array}{c}
\alpha(q-1)+\beta(q+1) = w\\
\beta \equiv m \pmod{q-1}.
\end{array}\right.
\end{equation*}
It is now a simple exercise to solve this system of
equations for each of the cases (i) to (ix). We leave it to
the reader. \CVD

\noindent\emph{Proof of Theorem~\ref{theo_calculs}.}
(1) We first prove the formulas involving the function $E$, and we
begin by showing (i). Let $n$ be an integer such that $0\le n< q$,
and suppose first that $n$ has the form $n=p^i$.
If $p\le p^i<q$ and $(p,q,i)\not=(2,4,1)$,
Proposition~\ref{derivees_mod_modulaires} (i) and
Lemma~\ref{espacesM} (i) yield immediately
$D_{p^i}E=E^{p^i+1}$. If $(p,q,i)=(2,4,1)$, the same Proposition and
the same Lemma now imply the existence of $\lambda\in C$ such that
$D_{p^i}E-E^{p^i+1}=\lambda g^{2}$. But since $D_{p^i}E-E^{p^i+1}$
vanishes at infinity (Lemma~\ref{lemme7}),
whereas $g$ does not, we have $\lambda=0$ and $D_{p^i}E=E^{p^i+1}$
again. Thus, the formula (i) is proved when $n$ has the form $n=p^i$
(for $i=0$ this follows from (\ref{system})). It is easy to deduce  
from this, by induction on $\alpha$, that the formula (i) also holds
when $n=\alpha p^i$ with $0\le\alpha\le p-1$ and $1\le p^i<q$. For if
the formula is true for $(\alpha-1)p^i$ (where $\alpha>1)$, then
\begin{align*}
\alpha D_{\alpha p^i}E=
D_{p^i}\circ D_{(\alpha-1)p^i}E & = D_{p^i}(E^{(\alpha-1)p^i}\cdot  E)\\
& = E^{(\alpha-1)p^i}D_{p^i}E+D_{p^i}(E^{(\alpha-1)p^i})E\\
& = E^{\alpha p^i+1}+(D_1E^{\alpha-1})^{p^i}E \\
& = E^{\alpha p^i+1}+ ((\alpha-1)E^{\alpha-2}(D_1E))^{p^i}E\\
& = \alpha E^{\alpha p^i+1}.
\end{align*}\label{pasdidee}
Now, if $n$ is arbitrary in the range $1\le n<q$, we argue
by induction on $n$ as follows. Define
$i\ge 0$ as the integer such that $p^i\le n<p^{i+1}$, and write
$n=\alpha p^i+\beta$, where $1\le\alpha\le p-1$ and $0\le\beta<p^i$.
By induction hypothesis we may assume that $D_{\beta}E=E^{\beta+1}$,
from which we deduce
$$ D_nE=D_{\beta}(D_{\alpha p^i}E)=D_{\beta}(E^{\alpha p^i+1})=
E^{\alpha p^i}D_{\beta}E=E^{\alpha p^i+\beta+1}=E^{n+1}.$$

To prove (iv), we use Proposition~\ref{derivees_mod_modulaires} (i) and
Lemma~\ref{espacesM} (ii).
In the case
$$ (p,q,i)\not\in\{(2,4,3),(3,3,1), (2,2,1)\},$$
we find that there exists $\lambda\in C$
such that $D_{p^i}E-E^{p^i+1}=\lambda g^{p^{i-e}-1}h^{p^{i-e}+1}$.
Since by Lemma~\ref{developp}
$$ D_{p^i}E=\frac{1}{d_1^{p^{i-e}+1}}t^{p^{i-e}+1}+ \cdots$$
and $g=1+\cdots$, $h=-1+\cdots$,
it follows that $\lambda = 1/d_1^{p^{i-e}}$, which gives (iv).
If $(p,q,i)=(2,4,3)$, the same method shows that there exist
$\lambda, \mu\in C$
such that $D_{p^i}E-E^{p^i+1}=\lambda g^6 +\mu gh^3$. The comparison
of the $t$-expansions in both sides implies $\mu=0$ and
$\lambda=1/d_1^{p^{i-e}}$, from which (iv) follows again. The other cases
are similar.

To establish (vii), we argue exactly in the same way, but using now
Lemmata~\ref{espacesM} (iii) and \ref{developp} (xi).

(2) Let us now prove the formulas involving the function $g$.
The formulas (ii) have already been proved
(see (\ref{system}) and Lemma~\ref{fonction_g}), so we have to prove
(v) and (viii). The same method as the one used in the proof of
(iv) above applies: We simply use Part (ii) of
Proposition~\ref{derivees_mod_modulaires} instead of Part (i),
Lemma~\ref{developp} for the $t$-expansions, and Parts (iv), (v)
of Lemma~\ref{espacesM} to know a basis of the space $M_{2p^i+q-1,p^i}$.

(3) Finally, we prove the formulas (iii), (vi) and (ix).
Let $n$ be an integer such that $0\le n< q$,
and first suppose that $n$ has the form $n=p^i$.
If $p\le p^i<q$ and $(p,q,i)\not=(2,4,1)$, then
Proposition~\ref{derivees_mod_modulaires} (iii) and Lemma~\ref{espacesM}~(vi)
yield at once the formula (iii) for $n=p^i$. If $(p,q,i)=(2,4,1)$,
arguing similarly and using the fact that $D_{p^i}h-E^{p^i}h$ vanishes at
infinity, we obtain again (iii). Thus, (iii) is proved when $n=p^i$
(if $p^i=1$ it follows from (\ref{system})).
Arguing now as in the proof of (i), we first deduce from this that (iii)
remains valid when $n=\alpha p^i<q$ with $0\le\alpha\le p-1$,
and then, as before, that the formula is true for all $n$ with $1\le  
n<q$.

To prove (vi) and (ix), we apply Proposition~\ref{derivees_mod_modulaires}
again: The function
$$D_{p^i}h-E^{p^i}h-E^q(D_{p^i-q}h) $$
is thus a modular form belonging to $M_{2p^i+q+1,p^i+1}$ for all $i\ge e$.
But this function has a $t$-expansion which is clearly an $O(t^2)$
by Lemma~\ref{lemme7}. So, using Lemma~\ref{espacesM} (parts (vii)
and (viii)), we see that there exist elements $\lambda,\mu\in C$
such that
$$D_{p^i}h-E^{p^i}h-E^q(D_{p^i-q}h)=\lambda
g^{p^{i-e}}h^{p^{i-e}+1}+\mu g^{2q+1}h^2,$$
where $\mu=0$ if $p^i<q^2$. With the help of Lemma~\ref{developp},
we can compute the first terms of the $t$-expansions on both sides. We
find $\lambda=-1/d_1^{p^{i-e}}$ if $p^i<q^2$ and, if $p^i=q^2$,
$\lambda=-(\frac{d_1}{d_2}+\frac{1}{d_1^q})$, $\mu=-1/d_2$.
To complete the proof, it remains to express $D_{p^i-q}h$ as a
polynomial in $E$, $g$, $h$. If $p^i=q$, we have $D_{p^i-q}h=h$
and we get immediately (vi). If $p^i>q$, then $D_{p^i-q}h=O(t^2)$ and
therefore, by Lemma~\ref{espacesM} (ix) and Lemma~\ref{fonction_h},
there exists $\nu\in C$ such that
$D_{p^i-q}h=\nu g^{p^{i-e}-1}h^{p^{i-e}}$. The $t$-expansion
of $D_{p^i-q}h$ given by Lemma~\ref{developp} (x) now yields
$\nu=1/d_1^{p^{i-e}-1}$, thus proving (vi) and (ix) . \CVD

\section{Proofs of Theorem~\ref{theoreme_h} and
Corollary \ref{theoreme_hh}.\label{section:groups}}

We consider the subfields of $F$
(\footnote{By \cite[Theorem~27.2]{MH},
the iterative derivation $(D_n)_{n\ge 0}$ extends to an
iterative higher derivation on $F$ in a unique way.
We denote this extension by $(D_n)_{n\ge 0}$ again.})
\[F_k=\bigcap_{i=0}^k\mbox{Ker}(D_{p^i}),\]
with $k\in\NN$
and we set $F_{-1}=F$.
Clearly, $F_{k-1}\supset F_{k}$ for all $k\in\NN$
and by (\ref{eq:op2}), one sees that 
$D_{p^k}(fg)=(D_{p^k}f)g+f(D_{p^k}g)$ for $f\in F_{k-1},g\in F$.
In particular, the restriction of $D_{p^k}$ to
$F_{k-1}$ is a derivation and thus, if $k\geq 0$, $x,y\in F_{k-1}$ and
$y\not=0$, then:
\[D_{p^k}\left(\frac{x}{y}\right)=\frac{y(D_{p^k}x)-x(D_{p^k}y)}{y^2}\]
(see Proposition
2.2 of \cite{MP} for further general properties of these subfields).

We also need to work with two families of subsets of $F^\times$.
If $k\in\NN$, we write
$$
\Psi_k =\{f\in F^\times\mbox{ such that }(D_{p^j}f)/f\in\widetilde{M}
\mbox{ for all }
0\leq j\leq k\},$$
and we define 
$$G_{k}:=\Psi_{k}\cap F_{k-1}^\times.$$ We also write
$\Psi_{-1}=G_{-1}=F^\times$.
Clearly, $\Psi_{k-1}\supset\Psi_{k}$ for all $k\in\NN$ and
$G_k\subset\Psi_{k}$ for all $k\geq -1$.
Moreover, $G_{k-1}\supset F_{k-1}^\times\supset G_{k}$ for all $k\in\NN$.
In Lemma \ref{Lemma_gruppo} of Section \ref{section:preliminaries}
we will prove that the $G_k$'s and the $\Psi_k$'s are multiplicative groups.

\medskip

\noindent \emph{Plan of Section \ref{section:groups}.}
In order to describe the content of this part, it is worth
outlining the structure of proof of Theorem \ref{theoreme_h} we intend to
develop.

By Proposition \ref{hstable}, $h\in\cap_{k=0}^\infty \Psi_k$, hence
$C^\times h^\ZZ\subset\cap_{k=0}^\infty \Psi_k$. We also observe that
Theorem~\ref{theoreme_h} is equivalent to:
\[\bigcap_{k=0}^\infty \Psi_k=C^\times h^\ZZ.\]
We will prove that for all $k\geq 0$:
\begin{equation}\label{eq:psik}
\Psi_k=F^\times_kh^\ZZ.\end{equation}
This result, combined with the results of Section
\ref{section:texp} (where we study $t$-expansions in $F^\times$),
%and makes Theorem \ref{theoreme_h} reachable.
will easily yield Theorem~\ref{theoreme_h}.

The proof of the equalities (\ref{eq:psik}) is by induction on $k\geq -1$;
the case $k=-1$ being trivial,
let us assume that $\Psi_{k-1}=F^\times_{k-1}h^\ZZ$ and
let $f$ be a quasi-modular form in $\Psi_k$. We know by induction 
hypothesis that there exists $n\in\ZZ$ such that
$fh^{n}$ is in $F^\times_{k-1}$ so that to explain the idea behind the
induction step we may restrict our attention to $f\in G_k=\Psi_k
\cap F^\times_{k-1}$.

Since $h^s\in\cap_{k=0}^\infty\Psi_k$, for all $s\in\ZZ$ we have
$$D_{p^k}(fh^{-sp^k})=a_sfh^{-sp^k}$$ with $a_s$ a quasi-modular form of weight
$2p^k$ in $F^\times_{k-1}$. The depth of $a_s$ is in general $\leq p^k$.
In Section \ref{section3} we determine some properties of the associated 
polynomials $P_{D_nf}$; by means of them, we observe that since
$a_s\in F^\times_{k-1}$,
its depth is divisible by $p^k$. Hence, $l(a_s)\in\{0,p^k\}$.

In section \ref{section:further} we characterise the integers $s$ such that
$l(a_s)<p^k$; their set is non-empty. For those integers $s$, $a_s$ is a
modular form because then $l(a_s)=0$.

In Section \ref{section:modularF}, we observe that modular forms in
$F_{k-1}^\times$ are very easy to describe (contrarily to quasi-modular forms).
Thanks to these considerations, we see that if
$s$ is such that $a_s$ is modular (and we will find such integers), then
$a_s=0$ hence proving (\ref{eq:psik}) for the index $k$.
Finally, Sections \ref{section4} and \ref{section:non_principal} serve
to end the proofs of the Theorem and Corollary in question.

This sketch of proof will hopefully help the reader to access the present
part of the text, but it must be considered
carefully, as several technical problems arise while taking at the place 
of $f$ quasi-modular, a more general element of $F^\times$.

\subsection{Preliminaries.\label{section:preliminaries}}

We begin with the following Lemma.

\begin{Lemme}
For all $k\geq -1$, $\Psi_k$ and $G_k$ are multiplicative subgroups of
$F^\times 
$.\label{Lemma_gruppo}\end{Lemme}
\noindent \emph{Proof.} This is clear for $k=-1$.
Let us write \[\Theta_{n}=
\{f\in F^\times\mbox{ such that }(D_{j}f)/f\in\widetilde{M}
\mbox{ for all }0\leq j\leq n\},\]
so that $\Theta_{p^j}=\Psi_j$ for all $j\in\NN$. We will prove by induction
on $n$ that for all $n\geq 0$, $\Theta_n$ is a subgroup of
$F^\times$ (this will imply
the required property for the $\Psi_k$'s). It is obviously true for
$n=0$, so we suppose now that $n\geq 1$ and that $\Theta_j$ is a subgroup
of $F^\times$ for all $j\le n-1$. We will show that
if $x,y\in\Theta_n$, then $xy^{-1}\in\Theta_n$.

Let $x,y$ be elements of $\Theta_n$. For all
$j=0,\ldots,n$ we have $D_jx=\alpha_jx $ and
$D_jy=\beta_jy$, with $\alpha_j,\beta_j\in \widetilde{M}$. Now, since
$\Theta_j$ is a subgroup of $F^\times$ for all $j=0,\ldots,n-1$,
we have $xy^{-1}\in\Theta_j$, and therefore there
exists $\gamma_j\in \widetilde{M}$ with $D_j(xy^{-1})=\gamma_jxy^{-1}$.
Thus, it only remains to show that $D_n(xy^{-1})/(xy^{-1})
\in\widetilde{M}$. But we have
\begin{eqnarray*}
\alpha_n x=D_nx=D_n(xy^{-1}\cdot y)&=&yD_n(xy^{-1})+\sum_{j=0}^{n-1}
D_j(xy^{-1})D_{n-j}y\\
&=&yD_n(xy^{-1})+\sum_{j=0}^{n-1}\gamma_j\beta_{n-j}x,
\end{eqnarray*}
from which it follows that
$$D_n(xy^{-1})/(xy^{-1})=\alpha_n-\sum_{j=0}^{n-1}\gamma_j\beta_{n-j}
\in\widetilde{M}.$$
Hence, $xy^{-1}\in\Theta_n$ and $\Theta_n$ is a  
multiplicative subgroup of $F^\times$ for all $n$. The property for
$G_k$ follows because $G_k=\Psi_{k}\cap F_{k-1}^\times$ for all $k\ge 0$.\CVD

\subsection{$t$-expansions.\label{section:texp}}

We obviously have an embedding $F\subset C((t))$.
More generally, we have the following Lemma.
\begin{Lemme}
For all $k\geq 0$ we have an embedding:
\[F_{k-1}\subset C((t^{p^k})).\]
\label{lemme3fede}\end{Lemme}
\noindent \emph{Proof.}
We check that for all $k\geq 0$, each element $x$ of $F_{k-1}$
has a convergent $t$-expansion (in a neighborhood of $0$):
\[x=\sum_{m\geq m_0}c_mt^{p^km},\quad c_m\in C.\]
We argue by induction on $k$, the case $k=0$ being clear.
Let $k$ be a strictly positive integer and let $x$ be in $F_{k-1}$.  
We then have $x\in F_{k-2}$, and by induction hypothesis:
\begin{equation}\label{doremi}
x=\sum_{m\geq m_0}c_mt^{p^{k-1}m}.
\end{equation}
By (\ref{eq:d1p}) and since $D_1t=t^2$, we have
$D_{p^{k-1}}(t^{p^{k-1}m})=(D_1t^m)^{p^{k-1}}=mt^{p^{k-1}(m+1)}$.
Thus, we get %(note that we may differentiate term-by-term
%by \cite[Proposition~2.6]{US}):
\begin{equation*}
D_{p^{k-1}}x=\sum_{m\geq m_0}mc_mt^{(m+1)p^{k-1}},
\end{equation*}
and $D_{p^{k-1}}x=0$ if and only if $mc_m=0$ for all $m$. This implies
$p\mid m$ if $c_m\not=0$, and hence $x\in C((t^{p^k}))$ by (\ref{doremi}).\CVD

\begin{Lemme}
Let $f$ be a non-zero element of $F$, and let us suppose that 
there exists a $p$-adic integer $s_1+s_2p+\cdots+s_{k+1}p^k+\cdots$
and a sequence $(f_k)_{k\in\NN}\subset F^\times$ such that, for all
$k$, $f_k\in F_k$ and
\begin{equation}f=h^{s_1+s_2p+\cdots+s_{k+1}p^k}f_k.\label{pippo}
\end{equation}
Then there exists $n\in\ZZ$ and $c\in C^\times$ with
$f=ch^{n}$.\label{lemme:sssss}
\end{Lemme}

\medskip

\noindent \emph{Proof.}
Since $F_{k}\subset C((t^{p^{k+1}}))$ (Lemma~\ref{lemme3fede}), we can write
\[f_k=b_kt^{\el_kp^{k+1}}+(\mbox{ terms of higher degree in $t$ }), 
\quad k=-1,0,\ldots,\] with $b_k\in C^\times$
and $\el_k\in\ZZ$. By (\ref{pippo}), we have
\[f_k=f_{k-1}h^{-s_{k+1}p^k},\quad k\geq 0.\]
Since $\nu_\infty(h)=1$, we obtain that:
\[\el_kp^{k+1}=\el_{k-1}p^k-s_{k+1}p^k,\quad k\geq 0,\] from which
we deduce
\begin{equation}
\el_k=\frac{\el_{k-1}}{p}-\frac{s_{k+1}}{p},\quad k\geq 0.
\label{pippa}\end{equation}

First, let us suppose that $\nu_\infty(f)=0$. Then
$\el_{-1}=0$, and (\ref{pippa}) implies $\el_0=\el_1=\cdots=0$  
because
for all $k$, $\el_k$ and $s_k$ are integers, and $0\leq s_{k+1}/p<1$ by  
hypothesis.
This implies $s_1=s_2=\cdots=0$ and $f=f_0=f_1=f_2=\cdots$. Hence:
\begin{equation*}f\in\bigcap_{k=0}^\infty F_{k}=C,
\end{equation*}
which implies $f=b_{-1}=c\in C^\times$.

Let us now suppose that $\nu_\infty(f)\not=0$. Since $\nu_\infty(h)=1$,
if we set $m=-\nu_\infty(f)$ and $\tilde{f}=h^mf$, then we have
$\nu_\infty(\tilde{f})=0$.
Let us now observe that there exist two sequences of integers
$(a_i)_{i\geq 1}$ and $(m_i)_{i\geq 0}$ unique
with the property that for all $k\geq 0$, $a_k\in\{0,\ldots,p-1 
\}$ and
\[m+s_1+s_2p+\cdots+s_{k+1}p^k=a_1+a_2p+\cdots+a_{k+1}p^k+m_kp^{k+1}.\]
The expression on the right hand side is nothing but the $p$-adic expansion
of the $p$-adic integer $m+\sum_{k=0}^\infty s_{k+1}p^k$, truncated to the
$k$-th digit. Hence, for all $k\geq 0$:
\begin{eqnarray*}
\tilde{f}=h^mf&=&h^{m+s_1+s_2p+\cdots+s_{k+1}p^k}f_{k}\\
&=&h^{a_1+a_2p+\cdots+a_{k+1}p^k}(h^{m_kp^{k+1}}f_k)\\
&=&h^{a_1+a_2p+\cdots+a_{k+1}p^k}\tilde{f}_k,
\end{eqnarray*}
with $\tilde{f}_k:=h^{m_kp^{k+1}}f_k$.
Since $h^{m_kp^{k+1}}\in F_k^\times$ and $f_k\in F_k^\times$,
we get $\tilde{f}_k\in F_{k}^\times$ for all $k$. Now,
since $\nu_\infty(\tilde{f})=0$, we can use the first part of the proof
and we obtain $\tilde{f}=h^mf=c\in C^\times$, whence $f=ch^{-m}$.\CVD

\subsection{Computing depths of derivatives.\label{section3}}

\noindent {\bf Terminology.}
In the following, we will say that a function \[\Phi:\Gamma_K\times 
\Omega\rightarrow C\] is
an {\em error term of degree $\leq l$} if for all $\gamma=\sqm{a}{b} 
{c}{d}\in\Gamma_K$ and $z\in\Omega$:
\[\Phi(\gamma,z)=\sum_{j=0}^l\sum_{i\in\ZZ/(q-1)\ZZ}f_{i,j}(z)\det(\gamma)^{- 
i}X_{\gamma}(z)^j,\]
where the $f_{i,j}$'s are holomorphic functions on $\Omega$ which  
do not depend on $\gamma$.

To simplify matters, the notation ${\cal R}(l)$ is reserved to  
designate any function which is an error term of degree $\leq l$. Hence,
two different functions may be denoted by the same symbol.

Let us denote by $\widetilde{M}_{w,*}^{l}$ the subset of $\widetilde 
{M}$ whose elements are
linear combinations of quasi-modular forms $f^{(i)}\in \widetilde{M}_ 
{w,i}^{\leq l}$ ($i\in\ZZ/(q-1)\ZZ$) such that there exists $i$
with $0\not=f^{(i)}\in\widetilde {M}_{w,i}^l$;
notice that this is not a $C$-vector space. By Theorem \ref{theo1},
$\widetilde{M}_{w,*}^{l}$ is also the subset of $C[E,g,h]$ of
the polynomials $P$  which are isobaric of weight $w$, such that
$\deg_E(P)=l$. If $P$ is an element of $C[E,g,h]$, we will also write
$l(P)=\deg_E(P)$, and we will refer to it as to the depth of $P$.

\begin{Proposition}
Let $f$ be a quasi-modular form of weight $w$, depth $l$ and
type $m$ (satisfying (\ref{QMcondition})),
let $n\geq 1$ be an integer. Then, for $\gamma\in\Gamma_K$, we have:
\begin{small}
\begin{eqnarray*}({\cal
D}_nf)(\gamma(z))&=&(cz+d)^{w+2n}\det(\gamma)^{-m-n}\times\\ & &\left 
(\binomial{w-l+n-1}{n}f_{l}(z)
X_{\gamma}(z)^{l+n}+{\cal R}(l+n-1)\right).
\end{eqnarray*}\end{small}

Let $f$ be a non-zero element of $\widetilde{M}_{w,*}^{l}$. Then:
\begin{small}
\begin{eqnarray*}({\cal
D}_nf)(\gamma(z))&=&(cz+d)^{w+2n}(\det\gamma)^{-n}\times\\
& &\left(\binomial{w-l+n-1}{n}\Pi_\gamma(z)
X_{\gamma}(z)^{l+n}+{\cal R}(l+n-1)\right),
\end{eqnarray*}\end{small}
where for all $\gamma\in\Gamma_K$, $\Pi_\gamma$ is a non-zero  
holomorphic function $\Omega\rightarrow C$ which does not depend on $n$.

The depth of ${\cal
D}_nf$ is $<l+n$ if and only if
\begin{equation}\binomial{w-l+n-1}{n}\equiv 0\pmod{p}.\label 
{eq:binomial0}\end{equation}
\label{lemme:derivate}\end{Proposition}
\noindent \emph{Proof.}
The first formula of the Proposition is a simple application of  
Proposition \ref{propo}.
We now prove the second part of the Proposition and we begin by  
explaining how to construct $\Pi_\gamma$.
Let us consider a non-zero element $f$ of $\widetilde{M}_{w,*}^{l}$.
We can write \begin{equation}
f=\sum_{i\in\ZZ/(q-1)\ZZ}f^{(i)}
\label{eq:indef2}
\end{equation}
with $f^{(i)}$ quasi-modular of weight $w$, depth $l_i$ and type $i$.
For $\gamma\in\Gamma_K$ we have:
\[f^{(i)}(\gamma(z))=(cz+d)^w\det(\gamma)^{-i}
\sum_{j=0}^{l_i}f^{(i)}_{j}(z)X_{\gamma}(z)^j\] with
$f^{(i)}_{j}$ quasi-modular for all $i,j$ (Lemma~\ref{lemme3}).
Thus, we can write:
\begin{eqnarray*}
f(\gamma(z))&=&(cz+d)^w
X_{\gamma}(z)^l\sum_{i}{}^*\det(\gamma)^{-i}f^{(i)}_{l}(z)+\\
&  & (cz+d)^w\sum_{i=0}^{q-2}\det(\gamma)^{-i}
\sum_{j=0}^{\min\{l_i,l-1\}}f^{(i)}_j(z)X_{\gamma}(z)^j\\
&=&(cz+d)^w\left(X_{\gamma}(z)^l\sum_{i}{}^*\det(\gamma)^{-i}f^{(i)}_{l}(z)
+{\cal R}(l-1)\right),
\end{eqnarray*}
where the sum $\sum^*$ runs over the indexes $i \in\ZZ/(q-1)\ZZ$ such that
$l(f^{(i)}) =l(f)=l$ (by hypothesis, this sum is non-empty).
By setting
\begin{equation}
\Pi_\gamma(z)=\sum{}^*\det(\gamma)^{-i}f^{(i)}_{l}(z)\label{eq:indef3}
\end{equation}
we get, for all
$\gamma=\sqm{a}{b}{c}{d}\in\Gamma_K$:
\begin{equation}
f(\gamma(z))=(cz+d)^w\left(X_{\gamma}(z)^l\Pi_\gamma(z)+{\cal R}(l-1) 
\right).
\label{eq:indef}
\end{equation}
By Lemma \ref{lemme3}, if $f^{(i)}_{l}\not=0$, it is a modular form of
weight $w-2l$ and type $i-l$.
For all $\gamma$, the function $\Pi_\gamma:\Omega\rightarrow C$ is  
not identically zero because, as follows from Proposition \ref{prop1}
(see also Theorem 5.13 of \cite{Ge}), non-zero modular forms of
the same weight but with different types are $C$-linearly independent.

Let us go back to the expression (\ref{eq:indef2}) for $f$; we have  
the identity (\ref{eq:indef}) with $\Pi_\gamma(z)$ as in (\ref 
{eq:indef3}).

More explicitly we have, for $\gamma=\sqm{a}{b}{c}{d}\in\Gamma_K$ and  
for $i\in\ZZ/(q-1)\ZZ$:
\[f^{(i)}(\gamma(z))=(cz+d)^w\det(\gamma)^{-i}(X_{\gamma}(z)^{l_i}
f_{l_i}^{(i)}(z)+{\cal R}(l_i-1)).\]
After the first part of the Proposition:
\begin{eqnarray*}
\lefteqn{({\cal D}_nf^{(i)})(\gamma(z))}\\ & = &(cz+d)^{w+2n}
(\det(\gamma))^{-i-n}\times\\ & &\left(\binomial{w-l_i+n-1}{n}f_{l_i}^{(i)} 
(z)X_{\gamma}(z)^{l_i+n}+{\cal R}(l_i+n-1)\right),\end{eqnarray*} thus
\begin{small}\begin{eqnarray*}
\lefteqn{({\cal D}_nf)(\gamma(z))}\\
& = & (cz+d)^{w+2n}(\det\gamma)^{-n}\left(\binomial{w-l+n-1}{n}
X_{\gamma}(z)^{l+n} 
\sum_i^*(\det(\gamma))^{-i}f_{l}^{(i)}+{\cal R}(l+n-1)\right)\\
& = &(cz+d)^{w+2n}(\det\gamma)^{-n}
\left(\binomial{w-l+n-1}{n}X_{\gamma}(z)^{l+n}\Pi_\gamma(z)+{\cal R}(l+n-1)
\right),
\end{eqnarray*}\end{small}
and the second part of the Proposition follows.

In particular, since $\Pi_\gamma$ is not  
identically zero for all $\gamma$,
the depth of ${\cal D}_nf$ is $<l+n$ if and only if (\ref{eq:binomial0})
holds.\CVD

\subsection{Further properties of elements in $\widetilde{M}
\cap F_k^\times$.\label{section:further}}

\begin{Lemme} For all $k\geq -1$, if $f\in\widetilde{M}_{w,*}^{l}\cap  
F_k^\times$, then $w-l=\alpha p^{k+1}$ with $\alpha\in\NN$.
\label{lemme:hhh}\end{Lemme}
\noindent \emph{Proof.}
We prove the Lemma by induction on $k$. Let $f$ be as in the hypotheses
of the Lemma. For $k=-1$ the result is trivial, so let us  
now suppose that $k\geq 0$, and that we have already proved the Lemma
for the integer $k-1$. As $F_{k}\subset F_{k-1}$, we have
$f\in F_{k-1}$ and by induction hypothesis we get $w-l=\alpha_{k-1}p^{k}$
with $\alpha_{k-1}\in\NN$. Moreover, by condition (\ref{eq:binomial0}) of
Proposition \ref{lemme:derivate},
\[0\equiv \binomial{w-l+p^k-1}{p^k}\equiv \binomial{(\alpha_{k-1}+1) 
p^k-1}{p^k}\pmod{p}.\]

Let $\beta=\beta_0+\beta_1 p+\cdots+\beta_sp^s$ be
a rational integer, with $\beta_0,\beta_1,\ldots\in\{0,\ldots,p-1\}$.
Since \[(1+\beta)p^k-1=(p-1)\sum_{i=0}^{k-1}p^i+
p^k(\beta_0+\beta_1 p+\cdots+\beta_sp^s),\]
taking into account (\ref{eq:lucas}) we see that
$\binomial{(\beta+1)p^k-1}{p^k}\equiv \beta\pmod{p}$. In particular,
\begin{equation}\binomial{(\alpha_{k-1}+1)p^k-1}{p^k}\equiv \alpha_ 
{k-1}\pmod{p}.\label{ertre}
\end{equation}
We thus obtain $\alpha_{k-1}=\alpha_{k}p$ with $\alpha_k\in\NN$, and
$w-l=\alpha_k p^{k+1}$.\CVD

\begin{Lemme}
Let $k\geq 0$ be an integer and let us consider elements
$f_1\in\widetilde{M}_{w_1,*}^{l_1}\cap F_{k-1}^\times$ and $f_2\in\widetilde 
{M}_{w_2,*}^{l_2}\cap F_{k-1}^\times$. For all $s\in\NN$ we have
\[f:=f_2h^{sp^k}({\cal D}_{p^k}f_1)-f_1{\cal D}_{p^k}(f_2h^{sp^k})\in 
\bigoplus_{m\in\ZZ/(q-1)\ZZ}\widetilde{M}^{\leq l_1+l_2+p^k}_{w_1+w_2 
+sp^k(q+1)+2p^k,m},\]
and there exists a unique integer $s\in\{0,\ldots,p-1\}$ such that:
\begin{equation}
l(f)<p^k+l_1+l_2.
\label{grossa_condizione}
\end{equation}
\label{lemme:svitamento}\end{Lemme}
\noindent \emph{Proof.} The first property of the Lemma easily  
follows from Theorem \ref{theo1}.

Since $f_1,f_2\in F_{k-1}^\times$ by hypothesis,
we have $w_1-l_1=\beta_1p^k$ and $w_2-l_2=\beta_2p^k$ with integers $
\beta_1,\beta_2\geq 0$ (Lemma~\ref{lemme:hhh}).
For all $\gamma\in\Gamma_K$ we have, after
Proposition \ref{lemme:derivate} or identity (\ref{eq:indef}):
\[
f_j(\gamma(z))=(cz+d)^{w_j}\left(\Pi_{\gamma,j}(z)X_{\gamma}(z)^{l_j}+
{\cal R}(l_j-1)\right),\quad j=1,2,
\]
where $\Pi_{\gamma,j}$ ($j=1,2$) is defined as in the equality
(\ref{eq:indef3}).
Proposition \ref{lemme:derivate} implies, for $\gamma=\sqm{a}{b}{c}{d}
\in\Gamma_K$:
\begin{eqnarray*}
\lefteqn{({\cal D}_{p^k}f_1)(\gamma(z))= (cz+d)^{w_1+2p^k}\times}\\
& &\left(\binomial{w_1-l_1+p^k-1}{p^k}\Pi_{\gamma,1}X_{\gamma}(z)^{l_1 
+p^k}+{\cal R}(l_1+p^k-1)\right),\\
\lefteqn{({\cal D}_{p^k}(f_2h^{sp^k}))(\gamma(z)) =(cz+d)^{w_2+sp^k(q
+1)+2p^k}\times}\\ & &
\left(\binomial{w_2+sp^k(q+1)-l_2+p^k-1}{p^k}\Pi_{\gamma,2}
(\det\gamma)^{-sp^k}h^{sp^k}X_{\gamma}(z)^{l_2+p^k}+\right.\\
& &\left.{\cal R}(l_2+p^k-1)\right).
\end{eqnarray*}
Hence,
\begin{small}
\begin{eqnarray*}\lefteqn{f_2(\gamma(z))h(\gamma(z))^{sp^k}({\cal D}_ 
{p^k}f_1)(\gamma(z))=}\\
&=&(cz+d)^{w(f)}\Biggl[\binomial{w_1-l_1+p^k-1}{p^k}\Pi_{\gamma,1}
\Pi_{\gamma,2}(\det\gamma)^{-sp^k}h^{sp^k}X_{\gamma}(z)^{l_1+l_2+p^k}+\\
& &{\cal R}(l_1+l_2+p^k-1)\Biggr],\\
\lefteqn{f_1(\gamma(z))({\cal D}_{p^k}(f_2h^{sp^k}))(\gamma(z))=
(cz+d)^{w(f)}\Biggl[\binomial{w_2+sp^k(q+1)-l_2+p^k-1}{p^k}\times}\\ &
&\Pi_{\gamma,1}\Pi_{\gamma,2}(\det\gamma)^{-sp^k}
h^{sp^k}X_{\gamma}(z)^{l_1+l_2+p^k}+{\cal R}(l_1+l_2+p^k-1)\Biggr].
\end{eqnarray*} \end{small}
This implies that for $\gamma\in\Gamma_K$:
\begin{small}
\[f(\gamma(z))=(cz+d)^{w(f)}\left(A\Pi_{\gamma,1}\Pi_{\gamma,2}
(\det\gamma)^{-sp^k} h^{sp^k}X_{\gamma}(z)^{l_1+l_2+p^k}
+{\cal R}(l_1+l_2+p^k-1)\right),\]
\end{small}
where \begin{eqnarray*}A&:=&\binomial{w_1-l_1+p^k-1}{p^k}- 
\binomial{w_2-l_2+p^ks(q+1)+p^k-1}{p^k}\\
&=&\binomial{p^k(\beta_1+1)-1}{p^k}-\binomial{p^k(\beta_2+s(q+1)+1)-1} 
{p^k}\\
&\equiv & \beta_1-(\beta_2+s)\pmod{p},
\end{eqnarray*}
after (\ref{ertre}).
Since the product:
\[(cz+d)^{w(f)}\Pi_{\gamma,1}\Pi_{\gamma,2}h^{sp^k}X_{\gamma}(z)^{l_1 
+l_2+p^k}\]
is not identically zero over $\Gamma_K\times\Omega$,
we have that $l(f)<l_1+l_2+p^k$ if and only if $s 
\equiv \beta_1-\beta_2\pmod{p}$, and there exists only one such integer in 
$\{0,\ldots,p-1\}$.\CVD

\medskip

\noindent {\bf Remark.} In particular, by choosing $f_1=1$ in Lemma  
\ref{lemme:svitamento}, we find the following property.
Let $k\geq 0$ be an integer and let $f\in\widetilde{M}_{w,*}^{l}\cap  
F_{k-1}$ be a non-zero quasi-modular form. There exists a unique integer
$s\in\{0,\ldots,p-1\}$ such that
\begin{equation}
l({\cal D}_{p^k}(h^{p^ks}f))<p^k+l.\label{eq:hallas}\end{equation}

\subsection{Modular forms in $F_k$\label{section:modularF}.}

We need some information about the modular forms contained in $F_k$.
\begin{Lemme}
Let $f\in C[g,h]$ be a modular form. If $D_1f=D_pf=\cdots=D_{p^k}f=0$,
then there exists a modular form
$\tilde{f}$ such that $f=\tilde{f}^{p^{k+1}}$.
\label{lemme:modulairesp}\end{Lemme}
\noindent \emph{Proof.} Let us consider the case $k=0$ first. By assumption,
$f$ is an element of $M_{w,m}$ for some
$w\in\NN$ and $m\in\ZZ/(q-1)\ZZ$. Let $\partial:C[g,h]\rightarrow C[g,h]$
be the derivation defined by $\partial:=-h\frac{\partial}{\partial g}$.
It is easily seen that the restriction of $\partial$
to $M_{w,m}$ coincides with the operator
$\partial_w:M_{w,m}\rightarrow M_{w,m}$ (of weight $2$
and type $1$) defined on p. 687 of \cite{Ge}, that is, we have
$\partial f=D_1f-wEf$.
Since $\partial f$ is a modular form and $D_1f=0$, we have thus
$\partial f=-wEf=0$, so that $w\equiv 0\pmod{p}$.
On the other hand, we have
$\partial f=-h\partial f/\partial g=0$, so $f$ is of the form
\begin{equation}
f=\sum_{{\tiny\begin{array}{c}i,j\\ \mbox{such that}\\ j\equiv 0\pmod 
{p}\end{array}}}\lambda_{i,j}h^ig^j.\label{eq:67}\end{equation}
Since the weight of $h$ is $q+1$ and $w\equiv 0\pmod{p}$, we see that  
$\lambda_{i,j}\not=0$ in (\ref{eq:67}) implies $i,j\equiv 0\pmod{p}$,
and $f$ is the $p$-th power of another modular form.

We continue the proof of the Lemma by induction on $k$: let us  
consider a modular form $f\in C[g,h]$ such that
$D_1f=D_pf=\cdots=D_{p^k}f=0$, with $k>0$.
Since $D_1f=D_pf=\cdots=D_{p^{k-1}}f=0$, we see from the induction  
hypothesis that there exists $r\in C[g,h]$ such that
$f=r^{p^k}$. Now, from (\ref{eq:d1p}) immediately follows:
\[0=D_{p^k}f=D_{p^k}(r^{p^k})=(D_1r)^{p^k}\] and $r$ is a $p$-th  
power after the first part of the Lemma:
$r=\tilde{f}^{p}$; thus $f=\tilde{f}^{p^{k+1}}$.\CVD

\begin{Lemme}  Let us suppose that $q\not\in\{2,3\}$.
For all $k\geq 0$ there does not exist a non-zero element of $ 
\widetilde{M}_{2p^k,*}^l\cap F_{k-1}$ with $l<p^{k}$.
\label{123}
\end{Lemme}
\noindent \emph{Proof.}
The Lemma for $k=0$ is clear: there does not exist a non-zero modular  
form of weight $2$, because $q\geq 4$.

Let us suppose by contradiction that for $k>0$ there exists a 
non-zero isobaric element $f$ of $\widetilde{M}\cap F_{k-1}$
of weight $2p^k$ and depth $l<p^k$.
Since $f\in F_{k-1}$ we have, after Lemma \ref{lemme:hhh}, $p^k(2- 
\alpha)=l$ with $\alpha\in\NN$,
and this condition implies $\alpha=2$ and $l=0$. Thus, $f$ is a  
modular form in $F_{k-1}$. After Lemma
\ref{lemme:modulairesp}, it is
a $p^k$-th power of another modular form $\tilde{f}$, which must have  
weight $2$, and which is zero after the first part of the Lemma; this  
implies a contradiction.\CVD

\noindent {\bf Remark.} For $q=2,3$, Lemma \ref{123} is false: $g$ has  
weight $2$ for $q=3$ and $g^2$ has weight $2$ for $q=2$; Lemma \ref 
{123} is the only tool of the proof of Theorem~\ref{theoreme_h}  
which needs the hypothesis $q\not\in\{2,3\}$.

\subsection{Proof of Theorem~\ref{theoreme_h}.\label{section4}}

We need the next Lemma.
\begin{Lemme}
For all $k\geq 0$ the following property holds. Let us consider $f\in  
G_k$: there exists a unique $s\in\{0,\ldots,p-1\}$ such that
\[h^{-p^ks}f\in F_{k}^\times.\] In particular,
\[G_k=h^{p^k\ZZ}F_{k}^\times.\]\label{lemme678}
\end{Lemme}
\noindent \emph{Proof.}
Write $f=a/b$ with $a,b\in \widetilde{M}$ non-zero. We have
$b^{p^k}f=ab^{p^k-1}$ with $f\in F_{k-1}^\times$ (by hypothesis)
and $b^{p^k}\in F_{k-1}$ (by (\ref{eq:d1p})), so
$x:=ab^{p^k-1}\in \widetilde{M}\cap F_{k-1}$ and
$y:=b^{p^k}\in \widetilde{M}\cap F_{k-1}$.
Thus we have $f=x/y$ with
$x,y\in \widetilde{M}\cap F_{k-1}^{\times}$.

Now, since $D_1h/h\in\widetilde{M}$ and $D_1(h^{p^i})=0$ for $i\ge 1$, we
easily deduce from (\ref{eq:d1p}) that $h^{p^k}\in G_k$. Hence, since $G_k$
is a group by Lemma \ref{Lemma_gruppo},
we find that for all $s\in\ZZ$, there exists $a=a_s\in \widetilde{M}$
such that:
\[D_{p^k}(h^{-sp^k}f)=ah^{-sp^k}f.\]
We claim that $a$ is isobaric of weight $2p^k$. To prove this,
we have of course to consider only the case $a\not= 0$.
Since $D_{p^k}$ is a derivation on $F_{k-1}$
and $h^{sp^k},x,y\in \widetilde{M}\cap F_{k-1}$, we have:
\[h^{sp^k}y(D_{p^k}x)-xD_{p^k}(h^{sp^k}y)=axyh^{sp^k}.\]
Let us write:
\[x=\sum_{i=\alpha}^\beta x_i,\quad
y=\sum_{j=\gamma}^\delta y_j,\quad
a=\sum_{l=\mu}^\nu a_l,\]
with $x_i,y_j,a_l\in\widetilde{M}$ isobaric of weights $i,j,l$, and  
$x_\alpha,x_\beta,y_\gamma,y_\delta,a_\mu,a_\nu$ non-zero. We have:
\begin{eqnarray}
\lefteqn{\sum_{v=\alpha+\gamma}^{\beta+\delta}\sum_{i+j=v}(y_jh^{sp^k} 
(D_{p^k}x_i)-x_iD_{p^k}(y_jh^{sp^k}))=}\nonumber\\
&=&
h^{sp^k}\sum_{w=\alpha+\gamma+\mu}^{\beta+\delta+\nu}\sum_{i+j+l=w} 
a_lx_iy_j.\label{eq:somme}
\end{eqnarray}
If a non-vanishing isobaric summand
$y_jh^{sp^k}(D_{p^k}x_i)-x_iD_{p^k}(y_jh^{sp^k})$ on the left-hand side
of (\ref{eq:somme})
is non-zero, then its weight $r:=2p^k+sp^k(q+1)+i+j$ satisfies $2p^k 
+sp^k(q+1)+\alpha+\gamma\leq r\leq 2p^k+sp^k(q+1)+\beta+\delta$.
Since $h^{sp^k}a_\mu x_\alpha y_\gamma$ and
$h^{sp^k}a_\nu x_\beta y_\delta$ are two non-vanishing isobaric
summands in the right-hand
side of (\ref{eq:somme}), whose weights are respectively
$sp^k(q+1)+\mu+\alpha+\gamma\leq sp^k(q+1)+\nu+\beta+\delta$,
we deduce that there exist two non-zero isobaric summands in the sum  
on the left-hand side.
We find:\begin{eqnarray*}
2p^k+sp^k(q+1)+\alpha+\gamma &\leq & sp^k(q+1)+\mu+\alpha+\gamma\\& 
\leq &sp^k(q+1)+\nu+\beta+\delta\\&\leq&
2p^k+sp^k(q+1)+\beta+\delta,\end{eqnarray*} which implies $2p^k\leq  
\mu\leq \nu\leq 2p^k$, that is, $\mu=\nu=2p^k$.
Hence $a\in \widetilde{M}$ is isobaric of weight $2p^k$, as claimed.

From this we deduce:
\begin{equation}h^{sp^k}y_\delta(D_{p^k}x_\beta)-x_\beta
D_{p^k}(h^{sp^k}y_\delta)=ax_\beta y_\delta h^{sp^k}.
\label{eq:datrattare}\end{equation}
After Lemma \ref{lemme:svitamento}, there exists exactly one
$s\in\{0,\ldots,p-1\}$
such that the term on the left-hand side of (\ref{eq:datrattare}) has  
depth $<p^k+l(x_\beta)+l(y_\delta)$.
Thus we see that the depth of $ax_\beta y_\delta h^{sp^k}$ is
$<p^k+l(x_\beta)+l(y_\delta)$.
Since the depth of $h$ is zero, we deduce that $a$ has depth $<p^k$.

Since $a=D_{p^k}(h^{-p^ks}f)/(h^{-p^ks}f)$ and
$h^{-p^ks}f,D_{p^k}(h^{-p^ks}f)\in F_{k-1}$,
we have $a\in F_{k-1}$.
Thus, $a\in \widetilde{M}\cap F_{k-1}$ is isobaric of weight $2p^k$  
and depth $<p^k$.
After Lemma \ref{123}, $a=0$; this implies $D_{p^k}(h^{-p^ks}f)=0$  
and $h^{-p^ks}f\in F_k$.
We deduce the equality of multiplicative groups
$G_k=h^{p^k\ZZ}F_{k}^\times$ as well.\CVD

We can terminate the proof of Theorem~\ref{theoreme_h}.
Let $f\in \bigcap_n\Psi_n$: thus, $f\in G_{0}=\Psi_1$.
After Lemma \ref{lemme678},
there exists (a unique) $s_1\in\{0,\ldots,p-1\}$ such that $f_0:=
h^{-s_1}f\in F_0^\times$; we continue by induction on $k$.

Let us suppose that
\[f_{k-1}:=h^{-(s_1+s_2p+\cdots+s_{k}p^{k-1})}f\in F_{k-1}^\times,\]  
for $k>0$ and $s_1,\ldots,s_{k}\in\{0,\ldots,p-1\}$.
Proposition \ref{hstable} and Lemma \ref{Lemma_gruppo}
imply that \[h^{-(s_1+s_2p+\cdots+s_{k}
p^{k-1})}\in \Psi_{k}\] and by hypothesis,
$f\in \Psi_{k}$. Thus, $f_{k-1}\in\Psi_{k}\cap F_{k-1}^\times=G_k$.

After Lemma \ref{lemme678}, there exists a unique $s_{k+1}\in\{0, 
\ldots,p-1\}$ such that
$f_k:=h^{-(s_1+s_2p+\cdots+s_{k+1}p^{k})}f\in F_{k}^\times$.

Thanks to this inductive process, for all $k\geq 0$, we can construct  
a sequence of elements $f_k\in F_k^\times$ and
a $p$-adic integer $s_1+s_2p+\cdots+s_{k+1}p^{k}+\cdots$ such that  
for all $k$, the equality
(\ref{pippo}) is satisfied. Lemma \ref{lemme:sssss} implies that  
there exists
$n\in\ZZ$ and $c\in C^\times$, with
$f=ch^{n}$.\CVD

\subsection{Proof of Corollary~\ref{theoreme_hh}.
\label{section:non_principal}}

If $U,V$ are isobaric polynomials of $C[E,g,h]$
of weights $w(U),w(V)$, their Rankin bracket $[U,V]$ is defined by:
\[[U,V]=\displaystyle{w(U) U (D_1V)-w(V) V (D_1U)}.\]
In the following Lemma
we collect the properties of these brackets that we need.

\begin{Lemme} Let $M$ be an isobaric element of $C[E,g,h]$.
\begin{itemize}
\item The map $d_M:X\mapsto [X,M]$ defined on the set of isobaric
elements of $C[E,g,h]$ satisfies
$$ d_M=d_M(E)\frac{\partial}{\partial E}+d_M(g)\frac{\partial}{\partial g}
+d_M(h)\frac{\partial}{\partial h}.$$
\item If ${\cal I}$ is an ideal of $C[E,g,h]$ such that
$D_1{\cal I}\subset{\cal I}$, and if $X\in{\cal I}$ is isobaric, then  
$d_M(X)\in{\cal I}$.
\end{itemize}\label{lemme:proprietes_crochet}\end{Lemme}

The proof of this Lemma is easy and left to the reader.
For the following Lemma, we recall that an ideal is {\em isobaric} if  
it is generated by isobaric elements of $C[E,g,h]$.

\begin{Lemme}
Every isobaric non-principal prime ideal ${\cal P}$ such that
$D_1{\cal P}\subset{\cal P}$ contains $h$.
\label{lemme:une_variable}\end{Lemme}
\noindent \emph{ Proof.}  We closely follow the proof of Lemma 2.3 of  
\cite{PF}.
Let ${\cal P}$ be a non-principal isobaric prime ideal of $C[E,g,h]$.
Eliminating $E$, we see that ${\cal P}$ contains a non-zero isobaric element
of $C[g,h]$, say $f$. We choose $f$ so that its total degree $\deg f$
is minimal. Let us define $f_g=\frac{\partial f}{\partial g}$ and
$f_h=\frac{\partial f}{\partial h}$. By minimality of $\deg f$,
$f$ is not a $p$-th power in $C[g,h]$. It follows that
$f_g\not=0$ or $f_h\not=0$.

We first suppose that $f_g\not=0$. By minimality of $\deg f$, we have
$f_g\not\in{\cal P}$. Now, by Lemma \ref{lemme:proprietes_crochet}, we
have $d_h(f)\in{\cal P}$ and
\begin{equation*}
d_{h}(f)  =  f_gd_{h}(g)+f_hd_{h}(h) = f_gd_{h}(g).
\end{equation*}
Since ${\cal P}$ is prime, it follows that $d_{h}(g)\in{\cal P}$.
But by (\ref{system}), $d_h(g)=h^2$. Thus, we get $h\in{\cal P}$.

If $f_g=0$ then $f_h\not=0$. By using the same arguments as before,
we find that $f_h\notin{\cal P}$, and then $d_g(h)=-h^2\in{\cal P}$,
which implies $h\in{\cal P}$.\CVD

\begin{Lemme}
Let  ${\cal I}$ be an ideal of $C[E,g,h]$ and let us denote by
$\tilde{{\cal I}}$ the ideal
generated by the isobaric elements of ${\cal I}$.
Let us suppose that ${\cal I}$ is hyperdifferential. Then, we have  
the following properties
\begin{enumerate}
\item The ideal $\tilde{{\cal I}}$ is hyperdifferential.
\item If ${\cal I}$ is prime, then $\tilde{{\cal I}}$ is prime.
\item If ${\cal I}$ is prime and non-principal, then $\tilde{{\cal I}} 
\not=(0)$.
\end{enumerate}
\label{lemme:nuovo}\end{Lemme}
\noindent \emph{ Proof.}
The first Part easily follows
from the fact that the images of any isobaric polynomial by the operators $D_n$
is again isobaric (Theorem \ref{theo2}).
The proof of the Parts two and three closely follows the proof of Lemma 5.2
of \cite{PF};
we do not need to give more details here.\CVD

\noindent \emph{Proof of Corollary~\ref{theoreme_hh}.}
Let ${\cal P}$ be non-zero hyperdifferential prime ideal.
If ${\cal P}$ is principal, then Theorem~\ref{theoreme_h} gives
immediately $h\in{\cal P}$, so we may suppose that ${\cal P}$
is non-principal.
By Lemma \ref{lemme:nuovo}, ${\cal P}$ contains a non-zero isobaric
hyperdifferential prime
ideal $\tilde{{\cal P}}$. If $\tilde{\cal P}$ is principal,
then $\tilde{{\cal P}}=(h)$ by Theorem~\ref{theoreme_h}.
If $\tilde{{\cal P}}$ is not principal, then
Lemma \ref{lemme:une_variable} implies
$h\in\tilde{{\cal P}}\subset{\cal P}$. \CVD

\section{Classification of hyperdifferential ideals.\label{section:examples}}

The aim of this Section is to provide a full description of the 
hyperdifferential ideals of $C[E,g,h]$. Let us write:
\begin{eqnarray*}
{\cal P}_0&=&(E ,h )\\
{\cal P}_\infty&=&(g,h)\\
{\cal P}_d&=&(h,E^{q-1}-dg),\quad d\in C^\times.
\end{eqnarray*}
These ideals are all prime, and we have the following diagram
of inclusions (here $c$ and $d$ vary in $C^\times$):
\begin{equation}\begin{array}{rcccl}(E,g-c,h)&\leftarrow& {\cal P}_0 & \\
& \swarrow & &\nwarrow&\\
(E,g,h) & \leftarrow &{\cal P}_d &\leftarrow&(h)\\
&\nwarrow& &\swarrow& \\
& & {\cal P}_\infty & &\end{array}\label{diagram}\end{equation}
 
In the following Theorem, we prove (under the condition $q\not=2,3$) that
all the ideals of the diagram (\ref{diagram}) are the only
non-zero hyperdifferential prime ideals of $C[E,g,h]$.

\begin{Theorem}\label{stable_ideals}

Let us assume that $q\not=2,3$.

\begin{itemize}
\item[(i)] If ${\cal P}$ is a principal hyperdifferential non-zero
prime ideal of $C[E,g,h]$, then ${\cal P}=(h)$.
\item[(ii)] For all $d\in C\cup\{\infty\}$,
the ideal ${\cal P}_d$ is hyperdifferential.
If ${\cal P}$ is a non-zero hyperdifferential prime ideal of $C[E,g,h]$ of
height $2$, then there exists $d\in C\cup\{\infty\}$ such that
${\cal P}={\cal P}_d$.
\item[(iii)] For all $c\in C$, the ideal $(E,g-c,h)$ is hyperdifferential.
If ${\cal P}$ is a maximal ideal which also is hyperdifferential, then there
exists $c\in C$ such that ${\cal P}=(E,g-c,h)$.
\end{itemize}
\end{Theorem}

Before going on with the proof of this Theorem, we will need two Lemmata.
In the next Lemma, we use the notation $f_1\equiv f_2\pmod{h}$
for $f_1,f_2\in C[E,g,h]$, which means that $f_2-f_1\in (h)$.

\begin{Lemme}\label{lemma:munu} For $n,\mu,\nu\in \NN$ we have
\begin{equation*}
D_n(E^\mu g^\nu)\equiv\binomial{\mu+\nu(q-1)+n-1}{n}E^{\mu+n}g^\nu\pmod{h}.
\end{equation*}
\end{Lemme}

\noindent\emph{Proof.} We have to show that for every monomial $f$ in $E$
and $g$, and for any $n\in\NN$, we have
\begin{equation}\label{munu}
D_n(f)\equiv \binomial{w(f)-l(f)+n-1}{n}E^nf \pmod{h}.
\end{equation}
First of all, we notice that if this formula holds for $f_1$ and $f_2$
(and all $n$), then it also holds for the product $f_1f_2$ : this follows
from Leibniz formula (\ref{eq:op2}) and Lemma~\ref{binom} with $N=n$,
$W=w(f)-l(f)$ and $M=l(g)-w(g)$
(notice also that $(-1)^i\binomial{M}{n-i}=(-1)^N\binomial{N-M-i-1}{N-i}$).
Thus, it suffices to prove the lemma for $f=E$ and for $f=g$.

Next, we show that it suffices to prove the formula (\ref{munu}) when $n$ has
the form $n=p^k$, $k\ge 0$. Indeed, assuming the formula for $n=p^k$,
an easy induction on $m$ first shows that for all $m\in\NN$,
$$ D_{p^k}^mf\equiv\left[\prod_{1\le i\le m}\binomial{w(f)-l(f)+ip^k-1}{p^k}
\right]E^{mp^k}f\pmod{h}$$
(here we need the fact that $f_1\equiv f_2\pmod{h}$ implies
$D_nf_1\equiv D_nf_2\pmod{h}$). Since
\begin{align*}
\prod_{1\le i\le m}\binomial{w(f)-l(f)+ip^k-1}{p^k} & =
\binomial{w(f)-l(f)+mp^k-1}{mp^k}\prod_{1\le i\le m}\binomial{ip^k}{p^k}\\
& = m!\binomial{w(f)-l(f)+mp^k-1}{mp^k}
\end{align*}
by (\ref{eq:lucas}), we find, using $D_{mp^k}=\frac{1}{m!}D_{p^k}^m$
(for $0\leq m\leq p-1$), that the formula (\ref{munu}) also holds for
$n=mp^k$ with $0\le m\le p-1$.
If now $n$ is arbitrary, write $m=n_0+\cdots+n_sp^s$ in base $p$
($0\le n_i\le p-1$).
Then formula~(\ref{munu}) for $n_ip^i$ together with (\ref{eq:op3}) gives
$$ D_nf=\left[\prod_{1\le i\le s}
\binomial{w(f)-l(f)+n_0+\cdots+n_ip^i-1}{n_ip^i}
\right]E^{n_0+\cdots+n_sp^s}f.$$
But the expression in brackets is also equal to
$$ \binomial{w(f)-l(f)+n-1}{n}\prod_{1\le i\le s}
\binomial{n_0+\cdots+n_ip^i}{n_ip^i},$$
which is equal to $\binomial{w(f)-l(f)+n-1}{n}$ by (\ref{eq:lucas}). So
(\ref{munu}) holds for all $n$.

It follows from the above remarks that we only need to prove the Lemma
when $n=p^k$ ($k\ge 0$) and for the functions $f=E$ and $f=g$. But in this
case this follows from Proposition~\ref{derivees_mod_modulaires},
Parts (i) and (ii).  Indeed, by Lemma \ref{lemme7}, the differences
$D_{p^k}E-E^{p^k+1}$ and $D_{p^k}g-E^{p^k}g$ are modular 
form which vanish at infinity, hence multiples of $h$.\CVD
 
The next Lemma, easy, will also be needed in the Proof of Theorem
\ref{stable_ideals}.

\begin{Lemme}
Let $r,r'$ be two rational integers such that
for all $n\geq 0$:
\[\binomial{r+n-1}{n}\equiv\binomial{r'+n-1}{n}\pmod{p}.\] Then, $r=r'$.
\label{lemma:rr1}\end{Lemme}
\noindent\emph{Proof.} Let $X$ be an indeterminate.
Let us suppose by contradiction that two distinct integers $r,r'$
exist, such that for all $n$, the congruence of the Lemma holds.
Since
\[(1-X)^{-r}=\sum_{n\geq 0}\binomial{r+n-1}{n}X^n\] in $\FF_p[[X]]$, and
a similar  identity holds for $(1-X)^{-r'}$, we have that
$(1-X)^{-r}=(1-X)^{-r'}$ in $\FF_p[[X]]$, which means that
$(1-X)^s=1$ with $s=r'-r\not=0$; this is impossible.\CVD

\noindent\emph{Proof of Theorem \ref{stable_ideals}: Part (i).}
This follows from Proposition \ref{hstable} and Theorem~\ref{theoreme_h}.

\medskip 

\noindent\emph{Part (ii).}
We introduce a new graduation on the
ring $C[E,g]$, namely the one defined by assigning to $E$ the degree $1$ and
to $g$ the degree $q-1$. In what follows, the word ``homogeneous''
will then refer to this new graduation.
%In other words, we ask that \[P=\sum_{\mu,\nu}c_{\mu,\nu}E^\mu g^\nu,\]
%where the coefficients $c_{\mu,\nu}$ are
%in $C$ and the (finite) sum runs over the pairs $(\mu,\nu)$ such that
%$r=\mu+\nu(q-1)$.

Let $P$ be a homogeneous polynomial in $C[E,g]$ of weight $r$.
By Lemma \ref{lemma:munu} we see that for all $n$,
\[D_nP\equiv \binomial{r+n-1}{n}E^nP\pmod{h}.\]
Thus, for such a polynomial $P$, the ideal $(P,h)$ is hyperdifferential.
In particular, the property above is fulfilled for the polynomial $E$ with
$r=1$, and for the polynomials $E^{q-1}-dg$ ($d\in C^\times$) and $g$,
with $r=q-1$. Hence the ideals ${\cal P}_0,{\cal P}_d$ (with $d\in C^\times$)
and ${\cal P}_\infty$ are hyperdifferential.

Let us now consider a hyperdifferential prime ideal ${\cal P}$ of
height $2$. By Theorem \ref{theoreme_hh}, $h\in{\cal P}$. Thus, there exists
an irreducible polynomial $P\in C[E,g]$ such that ${\cal P}=(P,h)$; without
loss of generality, we can suppose that $P$ is not proportional to $E$
(otherwise, ${\cal P}=(E,h)={\cal P}_0$ and we are done).
By assumption, for all $n\geq 0$ there exists $A_n\in C[E,g]$ such that
\begin{equation*}
D_nP\equiv A_nP\pmod{h}.
\end{equation*}

Since the operator $D_n:C[E,g]\rightarrow C[E,g,h]/(h)=C[E,g]$
is homogeneous of degree $n$ by Lemma~\ref{lemma:munu}, it is easy
to see that $A_n$ must be a linear combination of
homogeneous terms of weight $\leq n$.

We can write:
\[P=\sum_{r=0}^sP_r,\]
where $P_s\not=0$ and $P_r\in C[E,g]$ is homogeneous of weight $r$ for all $r$.
We have, explicitly by Lemma~\ref{lemma:munu},
($n\geq 0$, $0\leq r\leq s$)
$$
D_nP_r\equiv\binomial{r+n-1}{n}E^nP_r\pmod{h}.
$$
Thus, $E^n$ divides $A_nP$, which implies that $A_n=\lambda_n E^n$ for some $\lambda_n\in C$.
We get:
$$\sum_{r=0}^s\binomial{r+n-1}{n}P_r=\lambda_n\sum_{r=0}^sP_r,$$
which implies that for all $r$ such that $P_r\not=0$,
$$\binomial{s+n-1}{n}=\binomial{r+n-1}{n},\quad n\geq 0.$$

By Lemma \ref{lemma:rr1}, we find that the sum $\sum_{r=0}^sP_r$ has
exactly one non-vanishing term, $P_s$.
Therefore, $P$ is homogeneous of degree $s$, and we have:
\begin{equation}
P=E^\mu g^\nu\prod_{j=1}^a(E^{q-1}-\lambda_jg)^{\mu_j},\label{eq:product1}
\end{equation}
for integers $\mu,\nu,a,\mu_1,\ldots,\mu_a$ not all zero, and non-zero elements
$\lambda_1,\ldots,\lambda_a\in C$.
Now, $P$ is irreducible.
Hence, considering each type of factor in the product on the right hand side
of (\ref{eq:product1}), we get $P=E$ (case which is excluded) or $P=g$ or $P=E^{q-1}-dg$ for
$d\in C^\times$. If $P=g$ we get
${\cal P}={\cal P}_\infty$, and if $P=E^{q-1}-dg$ for some $d\in C^\times$,
then we get ${\cal P}={\cal P}_{d}$.

\medskip 

\noindent\emph{Part (iii).} Taking into account the $t$-expansions
(i), (ii) and (iii) of Lemma \ref{developp}, we first notice that
${\cal P}_0=(E,h)$ is also the ideal generated by the quasi-modular forms
which vanish at the infinity. After Lemma \ref{lemme7}, it contains the images
of all the operators $D_n$ ($n\geq 1$). Thus, we have $D_n(g-c)\in(E,h)$
for all $n\geq 1,c\in C$ (by the way, this argument gives an alternative
way to check that ${\cal P}_0$ is hyperdifferential).
Since now ${\cal P}_0$ is hyperdifferential after part (i),
it follows that for all $c\in C$, the maximal ideal $(E,g-c,h)$ is
hyperdifferential. 

Let ${\cal P}=(E-c_1,g-c_2,h-c_3)$ be a hyperdifferential maximal ideal.
In particular, $D_1E,D_1g,D_1h$ belong to ${\cal P}$, that is,
$E^2,-Eg-h,Eh\in{\cal P}$,
after the formulas (\ref{system}). Thus $E\in{\cal P}$, which implies $c_1=0$,
and $h\in{\cal P}$, which implies $c_3=0$.\CVD

\noindent\textbf{Remark.} There is also another way to check that
${\cal P}_\infty$ is hyperdifferential.
By Theorem \ref{theo1}, ${\cal P}_\infty$ is
the ideal generated by the quasi-modular forms $f$ which
have weights $w$ and depths $l$
satisfying the inequality $l<w/2$. By Theorem \ref{theo2}, this property
holds for $D_nf$, for all
$n\geq 1$ and for each quasi-modular form $f\in {\cal P}_\infty$. Hence,
$D_n {\cal P}_\infty\subset{\cal P}_\infty$ for all $n$.

\vspace{15pt}

\noindent{\small Vincent Bosser,\\
\indent Mathematisches Institut, Universit\"{a}t Basel,\\
\indent Rheinsprung 21,\\
\indent CH-4051 Basel, Switzerland.\\
\indent E-mail: {\tt Vincent.Bosser@unibas.ch}}

\vspace{15pt}

\noindent{\small Federico Pellarin,\\
\indent L.M.N.O., Universit\'e de Caen,\\
\indent Campus II - Boulevard Mar\'echal Juin,\\
\indent BP 5186 - F14032 Caen Cedex, France.\\
\indent E-mail: {\tt pellarin@math.unicaen.fr}}


\begin{thebibliography}{99}

\bibitem{abp} G.W. Anderson, W.D. Brownawell, M.A. Papanikolas.
{\em Determination of the algebraic relations among special $\Gamma$-values
in positive characteristic.} Ann. Math. 160, 237-313 (2004).

\bibitem{BV} V. Bosser. {\em Ind\'ependance alg\'ebrique de valeurs
de s\'eries d'Eisenstein}. In S.M.F. S\'eminaires et Congr\`es No. 12,
(2005).

\bibitem{changyu} C.-Y. Chang, J. Yu. {\em Determination of algebraic
relations among special zeta values in positive characteristic},
preprint (2006).

\bibitem{Dion} S. Dion.
{\em Un th\'eor\`eme de transcendance en caract\'eristique finie.}
J. Th\'eor. Nombres Bordx. 15, No. 1, 57-82 (2003).

\bibitem{KZ} M. Kaneko, D. Zagier.
{\em A generalized Jacobi theta function and quasimodular forms.}
Dijkgraaf, R. H. (ed.) et al., The moduli space of curves.  
%Proceedings of a conference held on Texel Island.
Basel: Birkh\"{a}user. Prog. Math. 129, 165-172 (1995).

\bibitem{Ge}  E.-U. Gekeler. {\em
On the coefficients of Drinfeld modular forms.}
Invent. Math. 93, No.3, 667-700 (1988).

%\bibitem{Ge1} E.-U. Gekeler.
%{\em Lectures on Drinfeld Modular Forms.}
%(Written by I. Longhi) CICMA Lecture Notes No. 4 (1999).

\bibitem{Lucas:Binomiaux} E. Lucas.
{\em Sur les congruences des nombres eul\'eriens et des
coefficients diff\'erentiels des fonctions trigonom\'etriques,
suivant un module premier.} Bull. Soc. Math. France, No. 6, pp. 49-54 (1878).

\bibitem{MH} H. Matsumura.
{\em Commutative ring theory.}
Cambridge Studies in Advanced Mathematics, 8. Cambridge etc.:  
Cambridge University Press. (1989).

\bibitem{MP} B. Matzat \& M. van der Put.
{\em Iterative differential equations and the Abhyankar conjecture.}
J. Reine Angew. Math. 557, 1-52 (2003).

\bibitem{MR} F. Martin \& E. Royer. {\em Formes modulaires et p\'eriodes.}
In S.M.F. S\'eminaires et Congr\`es No. 12, (2005).

\bibitem{Nes} Yu. V. Nesterenko. {\em Modular functions and transcendence
questions.} Sb. Math. 187, 1319-1348 (1996).

\bibitem{NP} Yu. V. Nesterenko et P. Philippon Editors.
{\em Introduction to algebraic independence theory.} Lecture Notes in  
Mathematics 1752, Springer (2001).

\bibitem{papa} M. A. Papanikolas. {\em Tannakian duality for Anderson-Drinfeld
motives and algebraic independence of Carlitz logarithms},
\texttt{http://arxiv.org/abs/math/0506078} (2005).

\bibitem{PF} F. Pellarin.
{\em La structure diff\'erentielle de l'anneau des formes 
quasi-modulaires pour ${\bf SL}_2({\mathbb Z})$.}
Journal de Th\'eorie des Nombres de Bordeaux 18, 241-264 (2006).

\bibitem{SD} H.P.F. Swinnerton-Dyer.
{\em On $\ell$-adic representations and congruences for coefficients of
modular forms.}
Modular Functions of one Variable III, Proc. internat. Summer School,
Univ. Antwerp 1972, Springer Lect. Notes Math. 350, 1-55 (1973).

\bibitem{US} Yu. Uchino, T. Satoh.
{\em Function field modular forms and higher derivations.}
Math. Ann. 311, No.3, 439-466 (1998).

\bibitem{vdp} J. Fresnel, M. Van der Put. {\em Rigid Analytic Geometry
and its Applications.} Progress in Mathematics 218, Birkh\"auser (2004).


\end{thebibliography}
\end{document}